\documentclass[11pt]{article}
\usepackage[english]{babel}

\usepackage{latexsym}
\usepackage{amsmath}
\usepackage{amsfonts}
\usepackage{mathrsfs}
\usepackage{amstext}
\usepackage{amssymb}
\usepackage{amsthm}       

\usepackage{xcolor}

\usepackage{hyperref}

\usepackage{tabu}
\usepackage{multirow}
\usepackage{array}

\usepackage{xspace}

\newtheorem{thm}{Theorem}[section]
\newtheorem{prop}[thm]{Proposition}

\newtheorem{lem}[thm]{Lemma}

\usepackage[margin=2cm]{geometry}

\title{The Dirichlet problem in a planar domain with two moderately close holes}
\author{ 
M.~Dalla Riva\thanks{Department of Mathematics, The University of Tulsa, 800 South Tucker Drive, Tulsa, Oklahoma 74104, USA.}$\ ^{,\dag}$ \   and P.~Musolino\thanks{Department of Mathematics, Aberystwyth University, Ceredigion SY23 3BZ, Wales, UK.}}

\date{\ }

\begin{document}
\maketitle

\noindent
{\bf Abstract:}  We investigate a Dirichlet problem for the Laplace equation in a domain of $\mathbb{R}^2$ with two small close holes. The domain is obtained by making in a bounded open set two perforations  at distance  $|\epsilon_1|$  one from the other and each one of size  $|\epsilon_1\epsilon_2|$.   In such a domain, we introduce a Dirichlet problem and we denote by $u_{\epsilon_1,\epsilon_2}$ its solution. We show that the dependence of $u_{\epsilon_1,\epsilon_2}$ upon $(\epsilon_1,\epsilon_2)$  can be described in terms of real analytic maps of the pair $(\epsilon_1,\epsilon_2)$  defined in an open neighborhood of $(0,0)$ and of logarithmic functions of $\epsilon_1$ and $\epsilon_2$.  Then we study the asymptotic behaviour of of $u_{\epsilon_1,\epsilon_2}$ as $\epsilon_1$ and $\epsilon_2$ tend to zero. We show that the first two terms of an asymptotic approximation can be computed only if we introduce a suitable relation between $\epsilon_1$ and $\epsilon_2$. 
\vspace{11pt}

 \par
\vspace{11pt}

%


\section{Introduction}

 The  asymptotic analysis of elliptic boundary value problems in domains with many  holes which collapse one to the other while shrinking their sizes is a topic of growing interest and  several authors have  recently  proposed  different techniques and points of view. We mention for example the method based on multiscale asymptotic expansions which have been used by    Bonnaillie-No\"el, Dambrine, Tordeux, and Vial \cite{BoDaToVi07, BoDaToVi09}, Bonnaillie-No\"el and Dambrine \cite{BoDa13}, and Bonnaillie-No\"el, Dambrine, and Lacave \cite{BoDaLa16}    to study problems with two moderately close holes, {\it i.e.}, problems with two holes whose mutual distance tends to zero while their size tends to zero at faster speed.  The  case when the number of holes is large has been considered by Maz'ya, Movchan, and Nieves  in a series of papers where they propose a  mesoscale approximation method to analyse problems for the Laplace operator and for the system of linear elasticity. We mention, for example, Maz'ya and Movchan \cite{MaMo10, MaMo12},  and Maz'ya, Movchan, and Nieves \cite{MaMoNi11, MaMoNi13, MaMoNi14, MaMoNi16}. The mesoscale approximation method does not require any periodicity assumption. If instead the holes have a periodic structure, then one can  resort to the large literature in homogenization  theory, where, rather then  aiming at obtaining asymptotic expansions, one typically characterizes the limit value of the solution of a perturbed problem as the solution of a limiting  problem. We refer, for instance, to the seminal works of Bakhvalov and Panasenko \cite{BaPa89}, Cioranescu and Murat \cite{CiMu82i,CiMu82ii}, and Mar{\v{c}}enko  and Khruslov \cite{MaKh74} and to the more recent `periodic unfolding method' used, {\it e.g.}, by Cioranescu, Damlamian, Donato,  Griso, and Zaki~\cite{CiDaDoGrZa12}).

In this paper, we consider  a Dirichlet problem for the Laplace equation in a planar domain with two small close holes. The method adopted   is different from those mentioned above. Indeed,  we follow the `functional analytic approach' which has been proposed by Lanza de Cristoforis for the analysis of linear and nonlinear singular perturbation problems (see, \textit{e.g.}, Lanza de Cristoforis \cite{La02,La10,La12}) and which allows the representation of the solution in terms of elementary functions and of real analytic maps of the singular perturbation parameters. One of the advantages of the method is that  real analytic maps   can be expanded into power series and thus, as a byproduct of our analysis, we can deduce fully justified asymptotic expansions for the solution with any order of approximation. Moreover, the coefficients of such expansions can be explicitly and constructively computed by solving certain systems of integral equations (as shown in \cite{DaMuRo15}). This method has been exploited for the analysis of Laplace and Poisson problems in domains with small close holes in \cite{DaMu16} and in \cite{DaMu}, respectively. In both of these papers, the conditions on the boundaries of the holes are of Neumann type. Here, instead, we will study a problem with Dirichlet conditions and we will focus on the two-dimensional case. This case is more involved than the higher dimensional case or the Neumann condition case because of the  logarithmic behaviour induced by the two-dimensional fundamental solution. As we shall see, such logarithmic behaviour will force the introduction of a specific relation between the size and the distance of the holes if we wish to  pass from the representation of the solution in terms of analytic maps to the explicit computation of the first asymptotic approximation terms.

We now proceed to introduce our problem and  we start by defining the geometric setting. We fix once for all a real number $\alpha\in]0,1[$ and three sets $\Omega^o$, $\Omega_1$ and $\Omega_ 2$ that satisfy the following condition: 
\begin{equation*}
\begin{split}
&\text{$\Omega^o$, $\Omega_1$ and $\Omega_ 2$ are open bounded connected subsets of $\mathbb{R}^2$}\\
&\text{of class $C^{1,\alpha}$, they contain the origin $0$ of $\mathbb{R}^2$ and}\\
&\text{they have connected boundaries $\partial\Omega^o$, $\partial\Omega_1$, and $\partial\Omega_2$.} 
\end{split}
\end{equation*}
Here the letter `$o$' stands for `outer domain' and $\Omega^o$ will play the role of the unperturbed outer domain in which we make two holes. To do so,  we take two points
 \[
 p^1 ,p^2 \in\mathbb{R}^2\,,\quad p^1 \neq p^2 
 \] 
 and  we assume that there exists 
\[
\delta_2>0
\]
such that 
\begin{equation}\label{empty1}
(p^1+\epsilon_2\, \mathrm{cl}\Omega_1)\cap(p^2+\epsilon_2\, \mathrm{cl}\Omega_2)=\emptyset\qquad\forall\epsilon_2\in[-\delta_2,\delta_2]\,.
\end{equation}
 Here and in the sequel `$\mathrm{cl}$' denotes the closure.   Then we define the rescaled sets
 \[
  \Omega_1(\epsilon_1,\epsilon_2)\equiv \epsilon_1 p^1 +\epsilon_1\epsilon_2 \Omega_ 1\,,\quad \Omega_2(\epsilon_1,\epsilon_2)\equiv \epsilon_1 p^2 +\epsilon_1\epsilon_2 \Omega_ 2\,,\quad \forall \epsilon_1,\epsilon_2\in\mathbb{R}\,,
 \]
 which will play the role of the holes. We observe that, for $\epsilon_1,\epsilon_2\in\mathbb{R}\setminus\{0\}$ and $i\in\{1,2\}$,  each $\Omega_i(\epsilon_1,\epsilon_2)$ is an open bounded subset of $\mathbb{R}^2$ which contains the point $\epsilon_1p^{i}$.  Instead, when $\epsilon_1=0$ or $\epsilon_2=0$, $\Omega_i(\epsilon_1,\epsilon_2)$ collapses to a point and we have  $\Omega_i(0,\epsilon_2)=\{0\}$ and  $\Omega_i(\epsilon_1,0)=\{\epsilon_1p^{i}\}$.  In addition, condition \eqref{empty1} implies that 
 \[
 \mathrm{cl}\Omega_ 1(\epsilon_1,\epsilon_2)\cap \mathrm{cl}\Omega_ 2(\epsilon_1,\epsilon_2) =\emptyset\qquad\forall\epsilon_1\in\mathbb{R}\setminus\{0\}\,,\;\epsilon_2\in[-\delta_2,\delta_2]\,.
 \]
 Then, one sees that   
 the mutual distance between $\Omega_1(\epsilon_1,\epsilon_2)$ and $\Omega_2(\epsilon_1,\epsilon_2)$ is controlled by $|\epsilon_1|$, while their size is proportional to $|\epsilon_1\epsilon_2|$. As a consequence, when both $\epsilon_1$ and $\epsilon_2$ approach zero, the size tends to zero at a faster rate than the distance. When this happens, one  says that the holes are `moderately close'.   In this paper,  we will also consider the case when the size and the distance are comparable, {\it i.e.} when $\epsilon_1$ tends to zero and $\epsilon_2$ stays away from zero. 
 
Since we want  the holes to be contained in $\Omega^o$,  we have to restrict the set of the `admissible' parameters $\epsilon_1$ for which we define the perforated domain. Then we take
\[
 \delta_1>0
 \]
such that 
\[
 \mathrm{cl}\Omega_1(\epsilon_1,\epsilon_2)\cup\mathrm{cl}\Omega_2(\epsilon_1,\epsilon_2)\subseteq\Omega^o\qquad\forall (\epsilon_1,\epsilon_2)\in[-\delta_1,\delta_1]\times[-\delta_2,\delta_2]
 \]
and we consider the pairs $(\epsilon_1,\epsilon_2)$ in the rectangular domain $[-\delta_1,\delta_1]\times[-\delta_2,\delta_2]$ as admissible parameters for which we define  the perforated domain
  \[
 \Omega(\epsilon_1,\epsilon_2)\equiv \Omega^o\setminus\left( \mathrm{cl}\Omega_1(\epsilon_1,\epsilon_2)\cup\mathrm{cl}\Omega_2(\epsilon_1,\epsilon_2)\right).
 \]
 We observe that for $\epsilon_1\in [-\delta_1,\delta_1]\setminus\{0\}$ and $\epsilon_2\in[-\delta_2,\delta_2]\setminus\{0\}$, $\Omega(\epsilon_1,\epsilon_2)$ is an open bounded connected subset of $\mathbb{R}^2$ of class $C^{1,\alpha}$ and the boundary of $\Omega(\epsilon_1,\epsilon_2)$ consists of three connected components: $\partial\Omega^o$, $\partial\Omega_1(\epsilon_1,\epsilon_2)$, and $\partial\Omega_2(\epsilon_1,\epsilon_2)$. For $\epsilon_1=0$ the set $\Omega(0,\epsilon_2)$ equals $\Omega^o\setminus\{0\}$ and for $\epsilon_2=0$ we have  $\Omega(\epsilon_1,0)=\Omega^o\setminus(\{\epsilon_1p^1 \}\cup\{\epsilon_1p^2 \})$. We also find convenient to introduce the notation
 \begin{equation}\label{tildaO}
 \tilde\Omega(\epsilon_2)\equiv\Omega_1(1,\epsilon_2)\cup\Omega_2(1,\epsilon_2)\qquad\forall\epsilon_2\in[-\delta_2,\delta_2]\,.
 \end{equation}

Now that the geometric configuration is settled, we turn to specify the boundary value problem. In order to define the Dirichlet data on $ \partial\Omega(\epsilon_1,\epsilon_2)$, we fix three functions
 \[
 f^o\in C^{1,\alpha}(\partial\Omega^o)\,,\; f_1\in C^{1,\alpha}(\partial\Omega_1)\,,\text{ and }f_2\in C^{1,\alpha}(\partial\Omega_2)\,.
\]
Then, for $\epsilon_1\in [-\delta_1,\delta_1]\setminus\{0\}$ and $\epsilon_2\in[-\delta_2,\delta_2]\setminus\{0\}$, we consider the following boundary value problem for a function $u\in C^{1,\alpha}(\mathrm{cl}\Omega(\epsilon_1,\epsilon_2))$:
 \begin{equation}\label{dir}
 \left\{
 \begin{array}{ll}
 \Delta u=0 &\text{in $\Omega(\epsilon_1,\epsilon_2),$}\\
 u=f^o &\text{on }\partial\Omega^o\,,\\
 u(x)=f_1((x-\epsilon_1p^1 )/(\epsilon_1\epsilon_2))&\forall x\in\partial\Omega_1(\epsilon_1,\epsilon_2)\,,\\
 u(x)=f_2((x-\epsilon_1p^2 )/(\epsilon_1\epsilon_2))&\forall x\in\partial\Omega_2(\epsilon_1,\epsilon_2)\,.
 \end{array}
 \right.
 \end{equation}
 As is well known the solution of problem \eqref{dir} exists and is unique. We denote such solution by $u_{\epsilon_1,\epsilon_2}$. Our aim is twofold:  first, we want to investigate the dependence of $u_{\epsilon_1,\epsilon_2}$ upon $\epsilon_1$ and $\epsilon_2$; then, we want to obtain asymptotic approximations of  of $u_{\epsilon_1,\epsilon_2}$ as $(\epsilon_1,\epsilon_2)$ tends to a degenerate value $(0,\gamma_0)$, with $\gamma_0\in[0,\delta_2[$. We will not consider the case when  $\epsilon_2$ tends to zero and $\epsilon_1$ tends to a non zero value,  which corresponds to the situation when the holes shrink to two distinct points. Such latter case has been largely investigated in literature (cf., {\it e.g.}, Maz'ya, Nazarov, and Plamenevskij \cite{MaNaPl00}). 
 
 Concerning the first  of the two goals, in Theorem \ref{ue1e2} we provide a representation (of suitable restrictions) of  $u_{\epsilon_1,\epsilon_2}$ and of the rescaled functions $u_{\epsilon_1,\epsilon_2}(\epsilon_1p^1+\epsilon_1\epsilon_2\,\cdot\,)$ and $u_{\epsilon_1,\epsilon_2}(\epsilon_1p^2+\epsilon_1\epsilon_2\,\cdot\,)$ in terms of real analytic functions of the pair $(\epsilon_1,\epsilon_2)$ and of the explicitly known functions $\log |\epsilon_1|$ and $\log|\epsilon_1\epsilon_2|$. The rescaled functions $u_{\epsilon_1,\epsilon_2}(\epsilon_1p^h+\epsilon_1\epsilon_2\,\cdot\,)$, with $h\in\{1,2\}$, describe the solution in proximity of the boundary of the holes and play an important role if one wants to compute quantities related to the solution, such as the energy integral. As a consequence of Theorem \ref{ue1e2} we see that, for $x\in\Omega^o\setminus\{0\}$ fixed and possibly shrinking $\delta_1$ and $\delta_2$, we have
\begin{equation}\label{anal}
u_{\epsilon_1,\epsilon_2}(x)=u^o(x)+\epsilon_1\epsilon_2\;U_x[\epsilon_1,\epsilon_2]+F[\epsilon_1,\epsilon_2]^t\,\Lambda(\epsilon_1,\epsilon_2)^{-1}\, V_x[\epsilon_1,\epsilon_2]
\end{equation}
for all $\epsilon_1\in]-\delta_1,\delta_1[\setminus\{0\}$ and  all $\epsilon_2\in]-\delta_2,\delta_2[\setminus\{0\}$, where $u^o$ is the solution of the unperturbed Dirichlet problem in $\Omega^o$ with boundary datum $f^o$, the functions $U_x$, $F$, and $V_x$ are real analytic from $]-\delta_1,\delta_1[\times]-\delta_2,\delta_2[$ to $\mathbb{R}$, and  $\Lambda(\epsilon_1,\epsilon_2)$ is a $2\times 2$  matrix such that
\[
\Lambda(\epsilon_1,\epsilon_2)\equiv 
\frac{1}{2\pi}\left(
\begin{array}{ll}
\log|\epsilon_1\epsilon_2|&\log|\epsilon_1|\\
\\
\log|\epsilon_1|&\log|\epsilon_1\epsilon_2|
\end{array}
\right)+R[\epsilon_1,\epsilon_2]
\]
 with $R$ real analytic from $]-\delta_1,\delta_1[\times]-\delta_2,\delta_2[$ to the space of $2\times 2$ real matrices. As we shall see, $\Lambda(\epsilon_1,\epsilon_2)$ is invertible if  both  $\epsilon_1$ and $\epsilon_2$ are not zero.

  Then, if we want to exploit \eqref{anal} to deduce asymptotic approximations of the solution as the pair $(\epsilon_1,\epsilon_2)$ approaches a degenerate value $(0,\gamma_0)$, we have to compute the inverse of the matrix $\Lambda(\epsilon_1,\epsilon_2)$. If we do so,  we obtain an expression which involves the quotient 
\begin{equation}\label{quotient}
\frac{\log |\epsilon_1|}{\log|\epsilon_1\epsilon_2|}
\end{equation}
(cf.~Proposition \ref{Lambda-1}). However, the limit of \eqref{quotient} as $(\epsilon_1,\epsilon_2)\to(0,\gamma_0)$ does not exist when $\gamma_0=0$.
 To overcome this difficulty, we introduce a relation between the parameters $\epsilon_1$ and $\epsilon_2$: we replace $\epsilon_1$ by a positive parameter $t$ and we take $\epsilon_2=\gamma(t)$, with $\gamma$ a function from a right neighbourhood of $0$ to $]0,\delta_2[$ such that the limits
\[
\gamma_0\equiv\lim_{t\to 0}\gamma(t)\qquad\text{and}\qquad\lambda_0\equiv\lim_{t\to 0}\frac{\log t}{\log(t\gamma(t))}
\]
exist finite in $[0,\delta_2[$ and $[0,+\infty[$, respectively. Under this assumption,  we obtain  in Proposition \ref{gamma} the first and second terms of  the asymptotic approximation of $u_{t,\gamma(t)}$ as $t>0$ tends to zero. In particular, for $\gamma_0=0$ we see that
\begin{equation}\label{asymp}
\begin{split}
&u_{t,\gamma(t)}(x)\\
&=u^o(x)+\frac{1}{\log(t\gamma(t))}\frac{2\pi}{1+\lambda_0}\Bigl(\lim_{y\to\infty} u_1(y)+\lim_{y\to\infty} u_2(y)-2u^o(0)\Bigr)G^{\Omega^o}(x,0)+o\left(\frac{1}{\log(t\gamma(t))}\right)
\end{split}
\end{equation}
as $t$ tends to zero. Here, $u_i$ with $i\in\{1,2\}$ denotes the harmonic solution of the exterior Dirichlet problem in $\mathbb{R}^2\setminus\Omega_i$ with boundary datum $f_i$ and $G^{\Omega^o}$ is the Green function of $\Omega^o$. We note that the limit value $\lambda_0$ appears explicitly in the second asymptotic terms in the right hand side of \eqref{asymp}.  In Proposition \ref{gamma} we also consider the case when $\gamma_0>0$ and the holes shrink their size and mutual distance at a comparable speed. In such a case we compute the expansion
\[
\begin{split}
&u_{t,\gamma(t)}(x)=u^o(x)\\
 & +\frac{2\pi}{\log t}\biggl(\lim_{y\to\infty} \tilde u(y)-u^o(0)+\biggl.\left(H^{2,1}_{\tilde\Omega(\gamma_0)}-H^{1,2}_{\tilde\Omega(\gamma_0)}\right)\int_{\Omega_1(1,\gamma_0)}\nu_{\Omega_1(1,\gamma_0)}\cdot\nabla\tilde u\, d\sigma\biggr) G^{\Omega^o}(\cdot,0)_{|\mathrm{cl}\Omega_M}+o\left(\frac{1}{\log t}\right)
\end{split}
\]
as $t$ tends to zero. Here $\tilde u$ is the harmonic solution of  a Dirichlet problem in the exterior domain $\mathbb{R}^2\setminus\tilde\Omega(\gamma_0)$ (see \eqref{tildeu}) and $H^{2,1}_{\tilde\Omega(\gamma_0)}$, $H^{1,2}_{\tilde\Omega(\gamma_0)}$ are quantities related to the Green function in the exterior domain  $\mathbb{R}^2\setminus\tilde\Omega(\gamma_0)$ (cf.~Proposition \ref{HxOe}). 

To conclude this introduction, we observe that our result justifies the introduction of specific relations between the size and the distance when dealing with the Dirichlet problem in a domain with moderately close small holes. Conditions of this type appear also in other papers on the topic. For example, in  \cite{BoDaLa16},  Bonnaillie-No\"el, Dambrine, and Lacave have considered a Poisson problem with Dirichlet conditions in a domain with two moderately close holes. To compute the asymptotic expansion of the solution,  they have assumed that   the distance behaves like the size to some power $\beta \in ]0,1[$. A condition which corresponds, with our notation, to the case when $\gamma(t)=t^{(1-\beta)/\beta}$ and the quotient  \eqref{quotient} is constant and equal to $1-\beta$. Another example can be found in \cite{MaMo10}, where Maz'ya and Movchan have analysed a Poisson problem with Dirichlet conditions in a domain with a large number small close holes. In such paper, it is assumed that  the size is smaller than the distance to the power $7/4$ (with our notation, $\gamma(t)<t^{3/4}$) in order to obtain uniform approximations of the solution and that the size is smaller than the square of the distance  (with our notation, $\gamma(t)<t$) to have approximations in $H^1$ norm (see also Maz'ya, Movchan, and Nieves \cite{MaMoNi13}).

\medskip
 
The present paper is organised as follows. In Section \ref{prel} we present some preliminary results on the solution of the Dirichlet problem in  a planar domain with many holes via potential theory. In Sections \ref{M1M2} and \ref{L} we  study some auxiliary integral operators that we use to convert problem \eqref{dir} into integral equations, while in Section \ref{H} we introduce some functions playing an important role in the description of the limiting behaviour of the solution $u_{\epsilon_1,\epsilon_2}$. In Section \ref{rep}, we prove Theorem \ref{ue1e2} on the representation of $u_{\epsilon_1,\epsilon_2}$ in terms of real analytic maps and known functions. In Section \ref{as}, we prove Proposition \ref{gamma}  where we analyse the asymptotic behaviour of the solution. 
   
 \section{The Dirichlet problem in a domain with many holes}\label{prel}
 In this section, we present some results of classical potential theory and we show how to exploit them in order to solve the Dirichlet problem for the Laplace equation in a domain with many holes. The construction of the solution that we present here will be then used to convert problem \eqref{dir} into equivalent integral equations. We start by denoting by $S$ the function from $\mathbb{R}^2\setminus\{0\}$ to $\mathbb{R}$ defined by
\[
S(x)\equiv
\frac{1}{2\pi}\log\,|x|
\qquad\forall x\in\mathbb{R}^2\setminus\{0\}\,.
\]
 As is well known, $S$ is a fundamental solution for the Laplace operator in $\mathbb{R}^2$.

Let $\mathcal{O}$ be an open bounded subset of $\mathbb{R}^2$ of class $C^{1,\alpha}$. Let $\phi\in C^{0,\alpha}(\partial\mathcal{O})$.  Then $v_\mathcal{O}[\phi]$ denotes the single layer potential with density $\phi$. Namely,
\[
v_\mathcal{O}[\phi](x)\equiv\int_{\partial\mathcal{O}}\phi(y)S(x-y)\,d\sigma_y\qquad\forall x\in\mathbb{R}^2\\,
\] where $d\sigma$ denotes the arc length element on $\partial\mathcal{O}$. As is well known, $v_\mathcal{O}[\phi]$ is a continuous function from $\mathbb{R}^2$ to $\mathbb{R}$ and  the restrictions $v^+_\mathcal{O}[\phi]\equiv v_\mathcal{O}[\phi]_{|\mathrm{cl}\mathcal{O}}$ and $v^-_\mathcal{O}[\phi]\equiv v_\mathcal{O}[\phi]_{|\mathbb{R}^n\setminus\mathcal{O}}$ belong to $C^{1,\alpha}(\mathrm{cl}\mathcal{O})$ and to $C^{1,\alpha}_{\mathrm{loc}}(\mathbb{R}^2\setminus\mathcal{O})$, respectively.  Here $C^{1,\alpha}_{\mathrm{loc}}(\mathbb{R}^2\setminus\mathcal{O})$ denotes the space of functions on $\mathbb{R}^2\setminus\mathcal{O}$ whose restrictions to $\mathrm{cl}\mathcal{B}$ belong to  $C^{1,\alpha}(\mathrm{cl}\mathcal{B})$ for all open bounded subsets $\mathcal{B}$ of $\mathbb{R}^2\setminus\mathcal{O}$.

If  $\psi\in C^{1,\alpha}(\partial\mathcal{O})$, then $w_\mathcal{O}[\psi]$ denotes the double layer potential with density $\psi$. Namely,
\[
w_\mathcal{O}[\psi](x)\equiv-\int_{\partial\mathcal{O}}\psi(y)\;\nu_{\mathcal{O}}(y)\cdot\nabla S_n(x-y)\,d\sigma_y\qquad\forall x\in\mathbb{R}^2\,,
\] where $\nu_\mathcal{O}$ denotes the outer unit normal to $\partial\mathcal{O}$  and the symbol `$\cdot$' denotes the scalar product in $\mathbb{R}^2$. The restriction $w_\mathcal{O}[\psi]_{|\mathcal{O}}$ extends to a function $w^+_\mathcal{O}[\psi]$ of $C^{1,\alpha}(\mathrm{cl}\mathcal{O})$ and  the restriction $w_\mathcal{O}[\psi]_{|\mathbb{R}^n\setminus\mathrm{cl}\mathcal{O}}$ extends to a function $w^-_\mathcal{O}[\psi]$ of $C^{1,\alpha}_{\mathrm{loc}}(\mathbb{R}^2\setminus\mathcal{O})$.

 Let 
\[
W_\mathcal{O}[\psi](x)\equiv -\int_{\partial\mathcal{O}}\psi(y)\;\nu_{\mathcal{O}}(y)\cdot\nabla S_n(x-y)\, d\sigma_y\qquad\forall x\in\partial\mathcal{O}\,,
\] for all $\psi\in C^{0,\alpha}(\partial\mathcal{O})$, and
\[
W^*_\mathcal{O}[\phi](x)\equiv \nu_{\mathcal{O}}(x)\cdot\int_{\partial\mathcal{O}}\phi(y)\;\nabla S_n(x-y)\, d\sigma_y\qquad\forall x\in\partial\mathcal{O}\,,
\] for all $\phi\in C^{1,\alpha}(\partial\mathcal{O})$. As is well known (cf.~Schauder \cite{Sc31, Sc32}) $W_\mathcal{O}$ is compact from   $C^{1,\alpha}(\partial\mathcal{O})$  to itself and $W^*_\mathcal{O}$  is compact from  $C^{0,\alpha}(\partial\mathcal{O})$ to itself. In addition $W_\mathcal{O}$ and $W^*_\mathcal{O}$ are adjoint with respect to the duality on $C^{1,\alpha}(\partial\mathcal{O})\times C^{0,\alpha}(\partial\mathcal{O})$ induced by $L^2({\partial\Omega})$  (cf.~Kress \cite{Kr99}). As a consequence,  one immediately deduces the validity of the following.

\begin{lem}\label{fredC}
The operators $\pm\frac{1}{2}I_\mathcal{O}+ W_\mathcal{O}$ are Fredholm of index $0$ from  $C^{1,\alpha}(\partial\mathcal{O})$ to itself. The operators $\pm\frac{1}{2}I_\mathcal{O}+ W^*_\mathcal{O}$ are Fredholm of index $0$ from  $C^{0,\alpha}(\partial\mathcal{O})$ to itself. The operator $\frac{1}{2}I_\mathcal{O}+ W^*_\mathcal{O}$ is the adjoint of $\frac{1}{2}I_\mathcal{O}+ W_\mathcal{O}$ and the operator  $-\frac{1}{2}I_\mathcal{O}+ W^*_\mathcal{O}$ is the adjoint of $-\frac{1}{2}I_\mathcal{O}+ W_\mathcal{O}$ with respect to the duality on $C^{1,\alpha}(\partial\mathcal{O})\times C^{0,\alpha}(\partial\mathcal{O})$ induced by $L^2({\partial\Omega})$.
\end{lem}

By exploiting the operators $W_\mathcal{O}$ and $W^*_\mathcal{O}$ we can write the jump formulas
\begin{equation}\label{jump}
w^\pm_\mathcal{O}[\psi]_{|\partial\mathcal{O}}=\pm\frac{1}{2}\psi+W_\mathcal{O}[\psi]\quad\text{ and }\quad\nu_\mathcal{O}\cdot\nabla v^\pm_\mathcal{O}[\phi]_{|\partial\mathcal{O}}=\mp\frac{1}{2}\phi+W^*_\mathcal{O}[\phi]
\end{equation}
which hold for all functions $\psi\in C^{1,
\alpha}(\partial\mathcal{O})$ and $\phi\in C^{0,
\alpha}(\partial\mathcal{O})$ (cf., {\it e.g.},  Folland \cite[Chap.~3]{Fo95}). If $\psi\in C^{1,\alpha}(\partial\mathcal{O})$ then we also have
\begin{equation}\label{nojump}
\nu_\mathcal{O}\cdot\nabla w^+_\mathcal{O}[\psi]_{|\partial\mathcal{O}}=\nu_\mathcal{O}\cdot\nabla w^-_\mathcal{O}[\psi]_{|\partial\mathcal{O}}\,.
\end{equation}

Now assume that $\mathcal{O}$ has $N$ connected components and $\mathbb{R}^2\setminus\mathrm{cl}\mathcal{O}$ has $K+1$ connected components and denote by $\mathcal{O}_1$, \dots, $\mathcal{O}_N$ the (bounded) connected components of $\mathcal{O}$ and by $\mathcal{O}^-_0$, $\mathcal{O}^-_1$, \dots, $\mathcal{O}^-_K$  the connected components of $\mathbb{R}^2\setminus\mathrm{cl}\mathcal{O}$. Since $\mathbb{R}^2\setminus\mathrm{cl}\mathcal{O}$ has a unique unbounded connected component we can assume that $\mathcal{O}^-_1$, \dots, $\mathcal{O}^-_K$ are bounded and that $\mathcal{O}^-_0$ is unbounded.  

In the sequel we exploit the following notation: if $\mathcal{X}$ is a subspace of $L^1(\partial\mathcal{O})$ then we denote by $\mathcal{X}_0$  the subspace of $\mathcal{X}$ consisting of the functions which have zero integral mean. 

Then we have the following classical lemma, where we describe the kernels of the integrals operator involved in the jump formulas in \eqref{jump} (cf., {\it e.g.},  Folland \cite[Chap.~3]{Fo95}).

\begin{lem}\label{ker}The following statements hold.
\begin{enumerate}
\item[(i)] The map from $\mathrm{Ker}(\frac{1}{2}I_\mathcal{O}+W^*_\mathcal{O})$ to $\mathrm{Ker}(\frac{1}{2}I_\mathcal{O}+W_\mathcal{O})$ which takes $\mu$ to $v[\partial\mathcal{O},\mu]_{|\partial\mathcal{O}}$ is bijective.
\item[(ii)] The map from $\mathrm{Ker}(-\frac{1}{2}I_\mathcal{O}+W^*_\mathcal{O})_0$ to $\mathrm{Ker}(-\frac{1}{2}I_\mathcal{O}+W_\mathcal{O})$ which takes $\mu$ to $v[\partial\mathcal{O},\mu]_{|\partial\mathcal{O}}$ is one to one. 
\item[(iii)]    $\mathrm{Ker}(\frac{1}{2}I_\mathcal{O}+W_\mathcal{O})$ consists of the functions from $\partial\mathcal{O}$ to $\mathbb{R}$ which are constant on $\partial\mathcal{O}^-_j$ for all $j\in\{1,\dots,K\}$ and which are identically equal to $0$ on $\partial\mathcal{O}^-_0$.
\item[(iv)] $\mathrm{Ker}(-\frac{1}{2}I_\mathcal{O}+W_\mathcal{O})$ consists of the functions from $\partial\mathcal{O}$ to $\mathbb{R}$ which are constant on $\partial\mathcal{O}_j$, for all $j\in\{1,\dots,N\}$. \item[(v)]  If $\phi\in\mathrm{Ker}(\frac{1}{2}I_\mathcal{O}+W^*_{\mathcal{O}})$ and $\int_{\partial\mathcal{O}_j^-}\phi\, d\sigma=0$ for all  $j\in\{1,\dots,K\}$, then $\phi=0$.
\item[(vi)]  If $\phi\in\mathrm{Ker}(-\frac{1}{2}I_\mathcal{O}+W^*_{\mathcal{O}})$  and $\int_{\partial\mathcal{O}_j}\phi\, d\sigma=0$ for all  $j\in\{1,\dots,N\}$, then $\phi=0$.
\item[(vii)] If $\phi\in\mathrm{Ker}(-\frac{1}{2}I_\mathcal{O}+W^*_{\mathcal{O}})_0$ and $v[\partial\mathcal{O},\phi]_{|\partial\mathcal{O}}$ is constant on $\partial\mathcal{O}$, then $\phi=0$.
\end{enumerate}
\end{lem}

Moreover, by Lemma \ref{ker} (i), (iii), and (v) we deduce the validity of the following.

\begin{lem}\label{tau} 
For each $i\in\{1,\dots,K\}$ there exists a unique function $\tau_i\in C^{0,\alpha}(\partial\mathcal{O})$ such that 
\[
(\frac{1}{2}I_\mathcal{O}+W^*_\mathcal{O})\tau_i=0\quad\text{and}\quad\int_{\partial\mathcal{O}^-_j}\tau_i d\sigma=\delta_{i,j}\quad\forall j\in\{1,\dots,K\}\,.
\]
The set $\{\tau_1,\dots,\tau_K\}$ is a basis for $\mathrm{Ker}(\frac{1}{2}I_\mathcal{O}+W^*_\mathcal{O})$ and the set $\{v_\mathcal{O}[\tau_1]_{|\partial\mathcal{O}},\dots,v_\mathcal{O}[\tau_K]_{|\partial\mathcal{O}}\}$ is a basis for $\mathrm{Ker}(\frac{1}{2}I_\mathcal{O}+W_\mathcal{O})$.
\end{lem}

In the sequel we denote by $\mathcal{X}_{\mathcal{O},i}$ the function from $\partial\mathcal{O}$ to $\mathbb{R}$ defined by 
\begin{equation}\label{X}
\mathcal{X}_{\mathcal{O},i}(x)\equiv\delta_{i,j}\qquad\forall i,j\in\{0,1,\dots,K\}\,,\; x\in\partial\mathcal{O}^-_j\,,
\end{equation}
where  $\delta_{i,j}$ is the Kronecker delta function.  By Lemma \ref{ker} (iii) it follows  that $\{\mathcal{X}_{\mathcal{O},1},\dots,\mathcal{X}_{\mathcal{O},K}\}$ is a basis for $\mathrm{Ker}(\frac{1}{2}I_\mathcal{O}+W_\mathcal{O})$. We also adopt the following notation, if $\Gamma$ is a one dimensional manifold in $\mathbb{R}^2$, then its one dimensional Lebesgue measure is denoted by
 $|\Gamma|$. Then we deduce the validity of the following. 

\begin{lem}\label{LambdaO}
Let $\Lambda_\mathcal{O}\equiv(\lambda^{i,j}_{\mathcal{O}})_{(i,j)\in\{1,\dots,K\}^2}$ be the real $K\times K$-matrix with entries $\lambda^{i,j}_{\mathcal{O}}$ defined by
\[
\lambda^{i,j}_{\mathcal{O}}\equiv\frac{1}{|\partial\mathcal{O}^-_j|}\int_{\partial\mathcal{O}}v_\mathcal{O}[\tau_i]\,\mathcal{X}_{\mathcal{O},j} d\sigma=\frac{1}{|\partial\mathcal{O}^-_j|}\int_{\partial\mathcal{O}^-_j}v_\mathcal{O}[\tau_i]\,d\sigma\quad\forall (i,j)\in\{1,\dots,K\}^2\,.
\] Then $\Lambda_\mathcal{O}$ is invertible and we have
$v_\mathcal{O}[\tau_i]_{|\partial\mathcal{O}}=\sum_{j=1}^K\lambda^{i,j}_{\mathcal{O}}\mathcal{X}_{\mathcal{O},j}$ for all $i\in\{1,\dots,K\}$.
\end{lem}

We are now ready to deduce the validity of the following Proposition \ref{ugO}, where we show how to construct the solution of the Dirichlet problem in a multiply perforated domain by solving some suitable integral equations.

\begin{prop}\label{ugO}
Let $g\in C^{1,\alpha}(\partial\mathcal{O})$.  Let $u\in C^{1,\alpha}(\mathrm{cl}\mathcal{O})$ be the unique function such that $\Delta u=0$ and $u_{|\partial\mathcal{O}}=g$. Then the following statements hold: 
\begin{enumerate}
\item[(i)] There exists and is unique a function $\mu\in C^{1,\alpha}(\partial\mathcal{O})$ such that 
\begin{equation}\label{ugO.eq1}
\left\{
\begin{array}{ll}
(\frac{1}{2}I_\mathcal{O}+W_\mathcal{O})\mu=g-\sum_{i=1}^K\Bigl(\int_{\partial\mathcal{O}}g\tau_i\,d\sigma\Bigr)\mathcal{X}_{\mathcal{O},i}\,,\\
\int_{\partial\mathcal{O}}\mu \mathcal{X}_{\mathcal{O},j}\,d\sigma=\int_{\partial\mathcal{O}^-_j}\mu\,d\sigma=0&\forall j\in\{1,\dots,K\};
\end{array}
\right.
\end{equation}
\item[(ii)] We have
\[
u(x)\equiv w_\mathcal{O}^+[\mu]+\sum_{i,j=1}^K\Bigl(\int_{\partial\mathcal{O}}g\tau_i\,d\sigma\Bigr)(\Lambda_\mathcal{O}^{-1})_{i,j}v_\mathcal{O}[\tau_j](x)\qquad\forall x\in\mathrm{cl}\mathcal{O}\,.
\]
\end{enumerate}
\end{prop}
\proof (i) By Lemma \ref{tau} one verifies that the right hand side of the first equation in \eqref{ugO.eq1} is orthogonal to 
$\mathrm{Ker}(\frac{1}{2}I_\mathcal{O}+W^*_\mathcal{O})$. Then the validity of the statement follows by Lemma \eqref{fredC} and by the standard properties of Fredholm operators. 

(ii) It is a consequence of  statement (i), of Lemma \ref{LambdaO}, of  \eqref{jump}, of the mapping properties of single and double layer potentials, and of the uniqueness of the solution of the Dirichlet problem. \qed

\section{The auxiliary maps $M_1$ and $M_2$}\label{M1M2}

Proposition \ref{ugO} shows how to construct the solution of the Dirichlet problem in two steps: first one constructs a basis for the kernel of the adjoint integral operator as in Lemma \ref{tau}, then one finds the solution of the system of integral equations of \eqref{ugO.eq1}. We want to exploit this approach for solving problem \eqref{dir}. Therefore, in this section, we perform the first of the two steps described above. Moreover, since our problem is defined in a domain which depends on $\epsilon_1$ and $\epsilon_2$, the integral equations delivered by Lemma \ref{tau} and Proposition \ref{ugO} will be defined on an $(\epsilon_1,\epsilon_2)$-dependent domain as well. As we are going to show, we will get rid of this dependence by performing a convenient change of variables.

We now introduce the auxiliary maps $M_1$ and $M_2$ representing the counterpart of Lemma \ref{tau}. For all $i\in\{1,2\}$  we denote by $M_i\equiv(M^o_i,M_{i,1}, M_{i,2}, M^c_i)$ the map from $]-\delta_1,\delta_1[\times{]-\delta_2,\delta_2[}\times C^{0,\alpha}(\partial\Omega^o)\times C^{0,\alpha}(\partial\Omega_1)\times C^{0,\alpha}(\partial\Omega_2)$ to $C^{0,\alpha}(\partial\Omega^o)\times C^{0,\alpha}(\partial\Omega_1)\times C^{0,\alpha}(\partial\Omega_2)\times\mathbb{R}^2$ defined by
\[
\begin{split}
&M^o_i[\epsilon_1,\epsilon_2,\rho^o_i,\rho_{i,1},\rho_{i,2}](x)\equiv [(\frac{1}{2}I_{\Omega^o}+W^*_{\Omega^o})\rho^o_ i](x)\\
& +\nu_{\Omega^o}(x)\cdot\sum_{h=1}^2\int_{\partial\Omega_ h}\nabla S(x-\epsilon_1p^h-\epsilon_1\epsilon_2\eta)\rho_{i,h}(\eta)\,d\sigma_\eta\qquad\forall x\in\partial\Omega^o\,,\\
&M_{i,h}[\epsilon_1,\epsilon_2,\rho^o_i,\rho_{i,1},\rho_{i,2}](\xi)\equiv [(-\frac{1}{2}I_{\Omega_h}+W^*_{\Omega_h})\rho_ {i,h}](\xi)\\
& +\epsilon_2\nu_{\Omega_h}(\xi)\cdot\int_{\partial\Omega_ k}\nabla S(p^h-p^k+\epsilon_2(\xi-\eta))\rho_{i,h}(\eta)\,d\sigma_\eta\\
& +\epsilon_1\epsilon_2\nu_{\Omega_h}(\xi)\cdot\int_{\partial\Omega^o}\nabla S(\epsilon_1p^h+\epsilon_1\epsilon_2\xi-y)\rho^o_{i}(y)\,d\sigma_y\qquad\forall h,k\in\{1,2\}\,,\;h\neq k\,,\; \xi\in\partial\Omega_h\,,\\
&M^c_i[\epsilon_1,\epsilon_2,\rho^o_i,\rho_{i,1},\rho_{i,2}](x)\equiv  \biggl(\int_{\partial\Omega_1}\rho_{i,1}\,d\sigma-\delta_{1,i}\,,\;\int_{\partial\Omega_2}\rho_{i,2}\,d\sigma-\delta_{i,2}\biggr)\,
\end{split}
\]
for all $(\epsilon_1,\epsilon_2,\rho^o_i,\rho_{i,1},\rho_{i,2})\in ]-\delta_1,\delta_1[\times]-\delta_2,\delta_2[\times C^{0,\alpha}(\partial\Omega^o)\times C^{0,\alpha}(\partial\Omega_1)\times C^{0,\alpha}(\partial\Omega_2)$.

Then, by a straightforward computation based on the rule of change of variable in integrals and by Lemma \ref{tau}, one deduces the validity of the following Proposition \ref{Mee=0}.

\begin{prop}\label{Mee=0}
Let $i\in\{1,2\}$. If $\epsilon_1\in ]-\delta_1,\delta_1[\setminus\{0\}$ and $\epsilon_2\in ]-\delta_2,\delta_2[\setminus\{0\}$, then we have 
\[
M_i[\epsilon_1,\epsilon_2,\rho^o_i,\rho_{i,1},\rho_{i,2}]=0
\]
if and only if
\[
\Bigl(\frac{1}{2}I_{\Omega(\epsilon_1,\epsilon_2)}+W^*_{\Omega(\epsilon_1,\epsilon_2)}\Bigr)\tau_i=0\quad\text{and}\quad \int_{\partial\Omega_j(\epsilon_1,\epsilon_2)}\tau_i\,d\sigma=\delta_{i,j}\quad\forall j\in\{1,2\}
\]
with $\tau_i\in C^{0,\alpha}(\partial\Omega(\epsilon_1,\epsilon_2))$ defined by
\begin{equation}\label{taui}
\tau_i(x)\equiv
\left\{
\begin{array}{ll}
\rho^o_i(x)&\forall x\in\partial\Omega^o\,,\\
\frac{1}{|\epsilon_1\epsilon_2|}\rho_{i,h}\Big(\frac{x-\epsilon_1p^h}{\epsilon_1\epsilon_2}\Big)&\forall h\in\{1,2\}\,,\; x\in\partial\Omega_h(\epsilon_1,\epsilon_2)\,.\\
\end{array}
\right.
\end{equation}
Moreover, there exists a unique triple $(\rho^o_i[\epsilon_1,\epsilon_2],\rho_{i,1}[\epsilon_1,\epsilon_2],\rho_{i,2}[\epsilon_1,\epsilon_2])\in C^{0,\alpha}(\partial\Omega^o)\times C^{0,\alpha}(\partial\Omega_1)\times C^{0,\alpha}(\partial\Omega_2)$ such that $M_i[\epsilon_1,\epsilon_2,\rho^o_i[\epsilon_1,\epsilon_2],\rho_{i,1}[\epsilon_1,\epsilon_2],\rho_{i,2}[\epsilon_1,\epsilon_2]]=0$.
\end{prop}

We now pass to consider the case when $\epsilon_2=0$ in Proposition \ref{Me0=0} and the case when $\epsilon_1=0$ and $\epsilon_2\neq 0$ in Proposition \ref{M0e=0}. The proofs of Propositions   \ref{Me0=0} and  \ref{M0e=0} can be effected by straightforward computations and by exploiting Lemma \ref{ker}.

\begin{prop}\label{Me0=0}
Let $i\in\{1,2\}$. If $\epsilon_1\in ]-\delta_1,\delta_1[$ (and $\epsilon_2=0$), then we have 
\[
M_i[\epsilon_1,0,\rho^o_i,\rho_{i,1},\rho_{i,2}]=0
\]
if and only if
\[
\left\{
\begin{array}{ll}
[(\frac{1}{2}I_{\Omega^o}+W^*_{\Omega^o})\rho^o_i](x)=-\nu_{\Omega^o}(x)\cdot\nabla S(x-\epsilon_1p^i)&\forall x\in\partial\Omega^o\,,\\
(-\frac{1}{2}I_{\Omega_h}+W^*_{\Omega_h})\rho_{i,h}=0\ \text{and}\  \int_{\partial\Omega_h}\rho_{i,h}\,d\sigma=\delta_{i,h}&\forall h\in\{1,2\}\,.
\end{array}
\right.
\]
Moreover, there exists a unique triple $(\rho^o_i[\epsilon_1,0],\rho_{i,1}[\epsilon_1,0],\rho_{i,2}[\epsilon_1,0])\in C^{0,\alpha}(\partial\Omega^o)\times C^{0,\alpha}(\partial\Omega_1)\times C^{0,\alpha}(\partial\Omega_2)$ such that $M_i[\epsilon_1,0,\rho^o_i[\epsilon_1,0],\rho_{i,1}[\epsilon_1,0],\rho_{i,2}[\epsilon_1,0]]=0$.
\end{prop}

We also observe that Proposition \ref{Me0=0} implies that 
\begin{equation}\label{rhoe0}
\rho_{i,j}[\epsilon_1,0]=0\qquad\forall i,j\in\{1,2\}\,,\;i\neq j
\end{equation} (see also Lemma \ref{ker} (vi)). In the following Proposition \ref{M0e=0} we exploit the definition of $\tilde\Omega(\epsilon_2)$ introduced in \eqref{tildaO} and we consider the case $\epsilon_1=0$.

\begin{prop}\label{M0e=0}
Let $i\in\{1,2\}$. If $\epsilon_2\in ]-\delta_2,\delta_2[\setminus\{0\}$ (and $\epsilon_1=0$), then we have 
\[
M_i[0,\epsilon_2,\rho^o_i,\rho_{i,1},\rho_{i,2}]=0
\]
if and only if
\[
\left\{
\begin{array}{ll}
[(\frac{1}{2}I_{\Omega^o}+W^*_{\Omega^o})\rho^o_i](x)=-\nu_{\Omega^o}(x)\cdot\nabla S(x)&\forall x\in\partial\Omega^o\,,\\
(-\frac{1}{2}I_{\tilde\Omega(\epsilon_2)}+W^*_{\tilde\Omega(\epsilon_2)})\tilde\rho_{i}=0\,,\\
 \int_{\partial\Omega_h(1,\epsilon_2)}\tilde\rho_{i}\,d\sigma=\delta_{i,h}&\forall h\in\{1,2\}\,,
\end{array}
\right.
\] 
with $\tilde\rho_i\in C^{0,\alpha}(\partial\tilde\Omega(\epsilon_2))$ defined by
\begin{equation}\label{tilderhoi}
\tilde\rho_i(x)\equiv\frac{1}{|\epsilon_2|}\rho_{i,h}\Bigl(\frac{x-p^h}{\epsilon_2}\Bigr)\qquad\forall h\in\{1,2\}\,,\; x\in\partial\Omega_h(1,\epsilon_2)\,. 
\end{equation}
Moreover, there exists a unique triple $(\rho^o_i[0,\epsilon_2],\rho_{i,1}[0,\epsilon_2],\rho_{i,2}[0,\epsilon_2])\in C^{0,\alpha}(\partial\Omega^o)\times C^{0,\alpha}(\partial\Omega_1)\times C^{0,\alpha}(\partial\Omega_2)$ such that $M_i[0,\epsilon_2,\rho^o_i[0,\epsilon_2],\rho_{i,1}[0,\epsilon_2],\rho_{i,2}[0,\epsilon_2]]=0$.
\end{prop}

Our aim is now to show that  $(\rho^o_i[\epsilon_1,\epsilon_2],\rho_{i,1}[\epsilon_1,\epsilon_2],\rho_{i,2}[\epsilon_1,\epsilon_2])$ depends analytically on $(\epsilon_1,\epsilon_2)$. In order to do so, we plan to apply the implicit function theorem for real analytic maps in Banach space. Thus, we need to show the real analyticity of  $M_i$ and the invertibility of the partial differential of $M_i$. We do that in the following technical Lemma \ref{dM}.

\begin{lem}\label{dM}
Let $i\in\{1,2\}$. The following statements hold.
\begin{enumerate}
\item[(i)] The map $M_i$ is real analytic from $]-\delta_1,\delta_1[\times]-\delta_2,\delta_2[\times C^{0,\alpha}(\partial\Omega^o)\times C^{0,\alpha}(\partial\Omega_1)\times C^{0,\alpha}(\partial\Omega_2)$ to $C^{0,\alpha}(\partial\Omega^o)\times C^{0,\alpha}(\partial\Omega_1)\times C^{0,\alpha}(\partial\Omega_2)\times\mathbb{R}^2$.
\item[(ii)] Let $(\bar\epsilon_1,\bar\epsilon_2,\bar\rho^o_i,\bar\rho_{i,1},\bar\rho_{i,2})\in ]-\delta_1,\delta_1[\times]-\delta_2,\delta_2[\times C^{0,\alpha}(\partial\Omega^o)\times C^{0,\alpha}(\partial\Omega_1)\times C^{0,\alpha}(\partial\Omega_2)$, then 
\begin{equation}\label{dM.eq1}
\partial_{(\rho^o_i,\rho_{i,1},\rho_{i,2})}M_i[\bar\epsilon_1,\bar\epsilon_2,\bar\rho^o_i,\bar\rho_{i,1},\bar\rho_{i,2}]
\end{equation}
(the partial differential of $M_i$ with respect to $(\rho^o_i,\rho_{i,1},\rho_{i,2})$ evaluated at $(\bar\epsilon_1,\bar\epsilon_2,\bar\rho^o_i,\bar\rho_{i,1},\bar\rho_{i,2})$) is an isomorphism from $C^{0,\alpha}(\partial\Omega^o)\times C^{0,\alpha}(\partial\Omega_1)\times C^{0,\alpha}(\partial\Omega_2)$ to $C^{0,\alpha}(\partial\Omega^o)\times C^{0,\alpha}(\partial\Omega_1)\times C^{0,\alpha}(\partial\Omega_2)\times\mathbb{R}^2$.
\end{enumerate}
\end{lem}
\proof
The validity of statement (i) follows by standard properties of integral  operators with real analytic kernels and with no singularity  (see, e.g., Lanza de Cristoforis and the second author \cite{LaMu13}) and by classical mapping properties of layer potentials (cf., {\it e.g.}, Miranda \cite{Mi65}).

To prove statement (ii) we observe that  the partial differential \eqref{dM.eq1}  is delivered by  
\[
\begin{split}
\partial_{(\rho^o_i,\rho_{i,1},\rho_{i,2})}M^o_i[\bar\epsilon_1,\bar\epsilon_2,\bar\rho^o_i,\bar\rho_{i,1},\bar\rho_{i,2}](\rho^o_i,\rho_{i,1},\rho_{i,2})&= M^o_i[\bar\epsilon_1,\bar\epsilon_2,\rho^o_i,\rho_{i,1},\rho_{i,2}]\,,\\
\partial_{(\rho^o_i,\rho_{i,1},\rho_{i,2})}M_{i,h}[\bar\epsilon_1,\bar\epsilon_2,\bar\rho^o_i,\bar\rho_{i,1},\bar\rho_{i,2}](\rho^o_i,\rho_{i,1},\rho_{i,2})&= M_{i,h}[\bar\epsilon_1,\bar\epsilon_2,\rho^o_i,\rho_{i,1},\rho_{i,2}]\quad\forall h\in\{1,2\}\,,\\
\partial_{(\rho^o_i,\rho_{i,1},\rho_{i,2})}M^c_{i}[\bar\epsilon_1,\bar\epsilon_2,\bar\rho^o_i,\bar\rho_{i,1},\bar\rho_{i,2}](\rho^o_i,\rho_{i,1},\rho_{i,2})&= \Big(\int_{\partial\Omega_1}\rho_{i,1}\, d\sigma\,,\; \int_{\partial\Omega_2}\rho_{i,2}\, d\sigma\Big)
\end{split}
\] 
for all $(\rho^o_i,\rho_{i,1},\rho_{i,2})\in C^{0,\alpha}(\partial\Omega^o)\times C^{0,\alpha}(\partial\Omega_1)\times C^{0,\alpha}(\partial\Omega_2)$.   By classical potential theory (cf.~Section \ref{prel}) and by a standard argument based on the theorem of change of variables in integrals one verifies that for all fixed $(g^o,g_1,g_2,c_1,c_2)\in C^{0,\alpha}(\partial\Omega^o)\times C^{0,\alpha}(\partial\Omega_1)\times C^{0,\alpha}(\partial\Omega_2)\times\mathbb{R}^2$ there exists and is unique a triple $(\rho^o_i,\rho_{i,1},\rho_{i,2})\in C^{0,\alpha}(\partial\Omega^o)\times C^{0,\alpha}(\partial\Omega_1)\times C^{0,\alpha}(\partial\Omega_2)$ such that 
\[
\partial_{(\rho^o_i,\rho_{i,1},\rho_{i,2})}M_i[\bar\epsilon_1,\bar\epsilon_2,\bar\rho^o_i,\bar\rho_{i,1},\bar\epsilon_{i,2}](\rho^o_i,\rho_{i,1},\rho_{i,2})=(g^o,g_1,g_2,c_1,c_2)\,.
\]
Then the validity of statement (ii) follows by the open mapping theorem.\qed

\vspace{\baselineskip}

Then, by a standard argument based on the implicit function theorem for real analytic maps (cf. Deimling \cite{De85}) we deduce the following Proposition \ref{rhoee}.

\begin{prop}\label{rhoee}
Let $i\in\{1,2\}$. Then the map from $]-\delta_1,\delta_1[\times]-\delta_2,\delta_2[$ to $C^{0,\alpha}(\partial\Omega^o)\times C^{0,\alpha}(\partial\Omega_1)\times C^{0,\alpha}(\partial\Omega_2)$ which takes $(\epsilon_1,\epsilon_2)$ to $(\rho^o_i[\epsilon_1,\epsilon_2],\rho_{i,1}[\epsilon_1,\epsilon_2],\rho_{i,2}[\epsilon_1,\epsilon_2])$ is real analytic. Moreover, the set of zeros of $M_i$ in $]-\delta_1,\delta_1[\times]-\delta_2,\delta_2[\times C^{0,\alpha}(\partial\Omega^o)\times C^{0,\alpha}(\partial\Omega_1)\times C^{0,\alpha}(\partial\Omega_2)$ coincides with the graph of $(\rho^o_i[\cdot,\cdot],\rho_{i,1}[\cdot,\cdot],\rho_{i,2}[\cdot,\cdot])$. 
\end{prop}

\section{The auxiliary map $L$}\label{L}

As we have done in the previous section for the counterpart of Lemma \ref{tau} for our problem \eqref{dir}, we now turn to consider the corresponding statement for the system in \eqref{ugO.eq1} of Proposition \ref{ugO}. Also in this case, we find convenient to perform a change of variables and to introduce the auxiliary map  $L\equiv(L^o,L_1, L_2)$ from $]-\delta_1,\delta_1[\times]-\delta_2,\delta_2[\times C^{1,\alpha}(\partial\Omega^o)\times C^{1,\alpha}(\partial\Omega_1)_0\times C^{1,\alpha}(\partial\Omega_2)_0$ to $C^{1,\alpha}(\partial\Omega^o)\times C^{1,\alpha}(\partial\Omega_1)\times C^{1,\alpha}(\partial\Omega_2)$ defined by
\[
\begin{split}
&L^o[\epsilon_1,\epsilon_2,\theta^o,\theta_1,\theta_2](x)\equiv [(\frac{1}{2}I_{\Omega^o}+W_{\Omega^o})\theta^o](x)\\
& +\epsilon_1\epsilon_2\sum_{i=1}^2\int_{\partial\Omega_i}\nu_{\Omega_i}(\eta)\cdot\nabla S(x-\epsilon_1p^i-\epsilon_1\epsilon_2\eta)\theta_{i}(\eta)\,d\sigma_\eta-f^o(x)\qquad\forall x\in\partial\Omega^o\,,\\
&L_h[\epsilon_1,\epsilon_2,\theta^o,\theta_1,\theta_2](\xi)\equiv [(-\frac{1}{2}I_{\Omega_h}+W_{\Omega_h})\theta_h](\xi)\\
&-w_{\Omega^o}[\theta^o](\epsilon_1p^h+\epsilon_1\epsilon_2\xi)\\
& -\epsilon_2\int_{\partial\Omega_ k}\nu_{\Omega_k}(\eta)\cdot\nabla S(p^h-p^k+\epsilon_2(\xi-\eta))\theta_k(\eta)\,d\sigma_\eta\\
& +f_h(\xi)-\int_{\partial\Omega^o}f^o\rho^o_h[\epsilon_1,\epsilon_2]\,d\sigma-\sum_{i=1}^2\int_{\partial\Omega_i}f_i\rho_{h,i}[\epsilon_1,\epsilon_2]\,d\sigma\\
&\qquad\qquad\qquad\qquad\qquad\qquad\qquad\qquad\qquad\qquad\forall h,k\in\{1,2\}\,,\;h\neq k\,,\; \xi\in\partial\Omega_h\,,
\end{split}
\]
for all $(\epsilon_1,\epsilon_2,\theta^o,\theta_1,\theta_2)\in ]-\delta_1,\delta_1[\times]-\delta_2,\delta_2[\times C^{1,\alpha}(\partial\Omega^o)\times C^{1,\alpha}(\partial\Omega_1)_0\times C^{1,\alpha}(\partial\Omega_2)_0$.
 
 Then, by a straightforward computation based on the rule of change of variable in integrals and by Proposition \ref{ugO}, one deduces the validity of the following Proposition \ref{Lee=0}.

 \begin{prop}\label{Lee=0}
If $\epsilon_1\in ]-\delta_1,\delta_1[\setminus\{0\}$ and $\epsilon_2\in]-\delta_2,\delta_2[\setminus\{0\}$, then we have 
\[
L[\epsilon_1,\epsilon_2,\theta^o,\theta_1,\theta_2]=0
\]
if and only if
\[
\left\{
\begin{array}{ll}
\left(\frac{1}{2}I_{\Omega(\epsilon_1,\epsilon_2)}+W_{\Omega(\epsilon_1,\epsilon_2)}\right)\phi=f-\sum_{k=1}^2\int_{\partial\Omega(\epsilon_1,\epsilon_2)}f\tau_k\,d\sigma\; \mathcal{X}_{\Omega(\epsilon_1,\epsilon_2),k}\\
\int_{\partial\Omega_j(\epsilon_1,\epsilon_2)}\phi\,d\sigma=0& \forall j\in\{1,2\}
\end{array}
\right.
\]
with $\phi, f\in C^{1,\alpha}(\partial\Omega(\epsilon_1,\epsilon_2))$ defined by
\[
\begin{split}
&\phi(x)\equiv
\left\{
\begin{array}{ll}
\theta^o(x)&\forall x\in\partial\Omega^o\,,\\
\theta_{h}\Big(\frac{x-\epsilon_1p^h}{\epsilon_1\epsilon_2}\Big)&\forall h\in\{1,2\}\,,\; x\in\partial\Omega_h(\epsilon_1,\epsilon_2)\,,\\
\end{array}
\right.\\
&f(x)\equiv
\left\{
\begin{array}{ll}
f^o(x)&\forall x\in\partial\Omega^o\,,\\
f_{h}\Big(\frac{x-\epsilon_1p^h}{\epsilon_1\epsilon_2}\Big)&\forall h\in\{1,2\}\,,\; x\in\partial\Omega_h(\epsilon_1,\epsilon_2)\,,\\
\end{array}
\right.
\end{split}
\]
and  $\mathcal{X}_{\Omega(\epsilon_1,\epsilon_2),k}$, $\tau_k$ defined as in \eqref{X} and  \eqref{taui}, respectively.

Moreover, there exists a unique triple $(\theta^o[\epsilon_1,\epsilon_2],\theta_1[\epsilon_1,\epsilon_2],\theta_2[\epsilon_1,\epsilon_2])\in C^{1,\alpha}(\partial\Omega^o)\times C^{1,\alpha}(\partial\Omega_1)_0\times C^{1,\alpha}(\partial\Omega_2)_0$ such that $L[\epsilon_1,\epsilon_2,\theta^o[\epsilon_1,\epsilon_2],\theta_1[\epsilon_1,\epsilon_2],\theta_2[\epsilon_1,\epsilon_2]]=0$.
\end{prop}

We now pass to consider the case when $\epsilon_2=0$ in Proposition \ref{Le0=0} and the case when $\epsilon_1=0$ and $\epsilon_2\neq 0$ in Proposition \ref{L0e=0}.

\begin{prop}\label{Le0=0}
If $\epsilon_1\in]-\delta_1,\delta_1[$ (and $\epsilon_2=0$), then we have
\[
L[\epsilon_1,0,\theta^o,\theta_1,\theta_2]=0
\]
if and only if 
\begin{equation}\label{Le0=0.eq1}
\left\{
\begin{array}{ll}
(\frac{1}{2}I_{\Omega^o}+W_{\Omega^o})\theta^o=f^o\,,\\
(-\frac{1}{2}I_{\Omega_h}+W_{\Omega_h})\theta_h=-f_h+\int_{\partial\Omega_h}f_h\,\rho_{h,h}[\epsilon_1,0]\,d\sigma&\forall h\in\{1,2\}\,.
\end{array}
\right.
\end{equation}

Moreover, there exists a unique triple $(\theta^o[\epsilon_1,0],\theta_1[\epsilon_1,0],\theta_2[\epsilon_1,0])\in C^{1,\alpha}(\partial\Omega^o)\times C^{1,\alpha}(\partial\Omega_1)_0\times C^{1,\alpha}(\partial\Omega_2)_0$ such that $L[\epsilon_1,0,\theta^o[\epsilon_1,0],\theta_1[\epsilon_1,0],\theta_2[\epsilon_1,0]]=0$.
\end{prop} 
\proof 
If $\theta^o$ satisfies the first equation of  the system \eqref{Le0=0.eq1}, then by the properties of adjoint operators, by Proposition \ref{Me0=0}, and by the definition of the double layer potential we have
\[
\begin{split}
&\int_{\partial\Omega^o} f^o\, \rho^o_h[\epsilon_1,0]\,d\sigma\\
&\quad=\int_{\partial\Omega^o} [(\frac{1}{2}I_{\Omega^o}+W_{\Omega^o})\theta^o]\, \rho^o_h[\epsilon_1,0]\,d\sigma\\
&\quad=\int_{\partial\Omega^o} \theta^o\, (\frac{1}{2}I_{\Omega^o}+W^*_{\Omega^o})\rho^o_h[\epsilon_1,0]\,d\sigma\\
&\quad=-\int_{\partial\Omega^o} \theta^o(y)\, \nu_{\Omega^o}(y)\cdot\nabla S(y-\epsilon_1p^h)\,d\sigma_y=-w_{\Omega^o}[\theta^o](\epsilon_1p^h)\,.
\end{split}
\]
Then the validity of the proposition follows by a straightforward computation based on the rule of change of variable in integrals, by equality \eqref{rhoe0}, and by a standard argument based on Lemma \ref{ker} and Proposition \ref{Me0=0}.
\qed
 
  \vspace{\baselineskip}
 
 Then we turn to consider the case $\epsilon_1=0$.
 
 \begin{prop}\label{L0e=0}
If $\epsilon_2\in]-\delta_2,\delta_2[$ (and $\epsilon_1=0$), then we have
\[
L[0,\epsilon_2,\theta^o,\theta_1,\theta_2]=0
\]
if and only if 
\begin{equation}\label{L0e=0.eq1}
\left\{
\begin{array}{ll}
(\frac{1}{2}I_{\Omega^o}+W_{\Omega^o})\theta^o=f^o\,,\\
(-\frac{1}{2}I_{\tilde\Omega(\epsilon_2)}+W_{\tilde\Omega(\epsilon_2)})\tilde\theta=-\tilde f+\sum_{h=1} ^2\left(\int_{\partial\tilde\Omega(\epsilon_2)}\tilde f\,\tilde\rho_h\,d\sigma\right)\mathcal{X}_{\tilde\Omega(\epsilon_2),h}\,,
\end{array}
\right.
\end{equation}
with $\tilde\theta, \tilde f\in C^{1,\alpha}(\partial\tilde\Omega(\epsilon_2))$ defined by
\[
\tilde\theta(x)\equiv
\theta_h\Big(\frac{x-p^h}{\epsilon_2}\Big)\,,\; \tilde{f}(x)\equiv
f_h\Big(\frac{x-p^h}{\epsilon_2}\Big)\,, \qquad\forall h\in\{1,2\}\,,\;  x\in\partial\Omega_h(1,\epsilon_2)\,,
\]
and $\mathcal{X}_{\tilde\Omega(\epsilon_2),h}$, $\tilde\rho_h$ defined as in  \eqref{X} and \eqref{tilderhoi}, respectively.

Moreover, there exists a unique triple $(\theta^o[0,\epsilon_2],\theta_1[0,\epsilon_2],\theta_2[0,\epsilon_2])\in C^{1,\alpha}(\partial\Omega^o)\times C^{1,\alpha}(\partial\Omega_1)_0\times C^{1,\alpha}(\partial\Omega_2)_0$ such that $L[0,\epsilon_2,\theta^o[0,\epsilon_2],\theta_1[0,\epsilon_2],\theta_2[0,\epsilon_2]]=0$.
\end{prop}
\proof 
If $\theta^o$ satisfies the first equation of  the system \eqref{L0e=0.eq1}, then by the properties of adjoint operators, by Proposition \ref{M0e=0}, and by the definition of the double layer potential we have
\[
\begin{split}
&\int_{\partial\Omega^o} f^o\, \rho^o_h[0,\epsilon_2]\,d\sigma\\
&\quad=\int_{\partial\Omega^o} [(\frac{1}{2}I_{\Omega^o}+W_{\Omega^o})\theta^o]\, \rho^o_h[0,\epsilon_2]\,d\sigma\\
&\quad=\int_{\partial\Omega^o} \theta^o\, (\frac{1}{2}I_{\Omega^o}+W^*_{\Omega^o})\rho^o_h[0,\epsilon_2]\,d\sigma\\
&\quad=-\int_{\partial\Omega^o} \theta^o(y)\, \nu_{\Omega^o}(y)\cdot\nabla S(y)\,d\sigma_y=-w_{\Omega^o}[\theta^o](0)\,.
\end{split}
\]
Then the validity of the proposition follows by a straightforward computation based on the rule of change of variable in integrals  and by a standard argument based on Lemma \ref{ker} and Proposition \ref{M0e=0}.
\qed

 In the following Proposition \ref{Lperp} we show an orthogonality property of the operator $L$.
 
 \begin{prop}\label{Lperp}
We have  
 \begin{equation}\label{Lrho=0}
 \begin{split}
& \int_{\partial\Omega^o} L^o[\epsilon_1,\epsilon_2,\theta^o,\theta_1,\theta_2]\,\rho^o_h[\epsilon_1,\epsilon_2]\,d\sigma\\
&\qquad-\sum_{k=1}^2\int_{\partial\Omega_k} L_k[\epsilon_1,\epsilon_2,\theta^o,\theta_1,\theta_2]\,\rho_{h,k}[\epsilon_1,\epsilon_2]\,d\sigma=0
 \end{split}
 \end{equation}
for all $h\in\{1,2\}$ and for all $(\epsilon_1,\epsilon_2,\theta^o,\theta_1,\theta_2)\in ]-\delta_1,\delta_1[\times]-\delta_2,\delta_2[\times C^{1,\alpha}(\partial\Omega^o)\times C^{1,\alpha}(\partial\Omega_1)_0\times C^{1,\alpha}(\partial\Omega_2)_0$.
 \end{prop}
 \proof
Let  $\epsilon_1\in ]-\delta_1,\delta_1[\setminus\{0\}$ and $\epsilon_2\in]-\delta_2,\delta_2[\setminus\{0\}$. Let $\phi, f\in C^{1,\alpha}(\partial\Omega(\epsilon_1,\epsilon_2))$ be defined as in Proposition \ref{Lee=0}. Then the validity of \eqref{Lrho=0} follows by equality 
 \[
 L[\epsilon_1,\epsilon_2,\theta^o,\theta_1,\theta_2]=(\frac{1}{2}I_{\Omega(\epsilon_1,\epsilon_2)}+W_{\Omega(\epsilon_1,\epsilon_2)})\phi-f+\sum_{k=1}^2\int_{\partial\Omega(\epsilon_1,\epsilon_2)}f\tau_k\,d\sigma\; \mathcal{X}_{\Omega(\epsilon_1,\epsilon_2),k}\,,
 \]
  by the orthogonality of $\mathrm{Ran}(\frac{1}{2}I_{\Omega(\epsilon_1,\epsilon_2)}+W_{\Omega(\epsilon_1,\epsilon_2)})$ and of $\mathrm{Ker}(\frac{1}{2}I_{\Omega(\epsilon_1,\epsilon_2)}+W^*_{\Omega(\epsilon_1,\epsilon_2)})$, by Proposition \ref{Mee=0}, and by a straightforward computation. 
 
 If at least one of $\epsilon_1$ and $\epsilon_2$ is $0$, then, by  Propositions  \ref{Le0=0} and \ref{L0e=0}, by the properties of adjoint operators, and by the definition of the double layer potential, we have 
 \begin{equation}\label{Lorhoo=0}
 \begin{split}
&\int_{\partial\Omega^o}L^o[\epsilon_1,\epsilon_2,\theta^o,\theta_1,\theta_2]\,\rho^o_h[\epsilon_1,\epsilon_2]\,d\sigma\\
&\quad= \int_{\partial\Omega^o}\Bigl[\Bigl(\frac{1}{2}I_{\Omega(\epsilon_1,\epsilon_2)}+W_{\Omega(\epsilon_1,\epsilon_2)}\Bigl)\theta^o\Bigr]\,\rho^o_h[\epsilon_1,\epsilon_2]\,d\sigma
-\int_{\partial\Omega^o}f^o\,\rho^o_h[\epsilon_1,\epsilon_2]\,d\sigma\\
&\quad= \int_{\partial\Omega^o}\theta^o\,\Bigl(\frac{1}{2}I_{\Omega(\epsilon_1,\epsilon_2)}+W^*_{\Omega(\epsilon_1,\epsilon_2)}\Bigr)\rho^o_h[\epsilon_1,\epsilon_2]\,d\sigma
-\int_{\partial\Omega^o}f^o\,\rho^o_h[\epsilon_1,\epsilon_2]\,d\sigma\\
&\quad= -\int_{\partial\Omega^o}\theta^o(y)\,\nu_{\Omega^o}(y)\cdot\nabla S(y-\epsilon_1p^h)\,d\sigma_y
-\int_{\partial\Omega^o}f^o\,\rho^o_h[\epsilon_1,\epsilon_2]\,d\sigma\\
&\quad=-w_{\Omega^o}[\theta^o](\epsilon_1p^h)-\int_{\partial\Omega^o}f^o\,\rho^o_h[\epsilon_1,\epsilon_2]\,d\sigma\,.
 \end{split}
 \end{equation}
If $\epsilon_2=0$, then the orthogonality of $\mathrm{Ran}(-\frac{1}{2}I_{\Omega_h}+W_{\Omega_h})$ and of $\mathrm{Ker}(-\frac{1}{2}I_{\Omega_h}+W^*_{\Omega_h})$ and equality 
\[
\int_{\partial\Omega_k}\rho_{h,k}[\epsilon_1,0]\,d\sigma=\delta_{h,k}\qquad\forall k\in\{1,2\}
\]
(cf. Proposition \ref{Me0=0}) imply that 
\begin{equation}\label{Lkrhoke0=0}
\int_{\partial\Omega_k}L_k[\epsilon_1,0,\theta^o,\theta_1,\theta_2]\,\rho_{h,k}[\epsilon_1,0]\,d\sigma=\delta_{h,k}\Bigl(-w_{\Omega^o}[\theta^o](\epsilon_1p^h)-\int_{\partial\Omega^o}f^o\,\rho^o_h[\epsilon_1,0]\,d\sigma\Bigr)\,.
\end{equation}
If instead $\epsilon_1=0$ and $\epsilon_2\neq 0$, then by the orthogonality of $\mathrm{Ran}(-\frac{1}{2}I_{\tilde\Omega(\epsilon_2)}+W_{\tilde\Omega(\epsilon_2)})$ and of $\mathrm{Ker}(-\frac{1}{2}I_{\tilde\Omega(\epsilon_2)}+W^*_{\tilde\Omega(\epsilon_2)})$ and equality 
\[
\int_{\partial\tilde\Omega(\epsilon_2)}\tilde\rho_{h}\,d\sigma=\delta_{h,k}\qquad\forall k\in\{1,2\}\,,
\]
where $\tilde\rho_h\in C^{0,\alpha}(\partial\tilde\Omega(\epsilon_2))$ is defined as in Proposition \ref{M0e=0}, we deduce that
\begin{equation}\label{Lkrhok0e=0}
 \begin{split}
&\sum_{k=1}^2\int_{\partial\Omega_k}L_k[0,\epsilon_2,\theta^o,\theta_1,\theta_2]\,\rho_{h,k}[0,\epsilon_2]\,d\sigma\\
&=\int_{\partial\tilde\Omega(\epsilon_2)}\Bigl[\Bigl(-\frac{1}{2}I_{\tilde\Omega(\epsilon_2)}+W_{\tilde\Omega(\epsilon_2)}\Bigr)\tilde\theta\Bigr]\tilde\rho_h\,d\sigma-w_{\Omega^o}[\theta^o](0)-\int_{\partial\Omega^o}f^o\,\rho^o_h[0,\epsilon_2]\,d\sigma\\
&=-w_{\Omega^o}[\theta^o](0)-\int_{\partial\Omega^o}f^o\,\rho^o_h[0,\epsilon_2]\,d\sigma\,.
 \end{split}
\end{equation}
Now the validity of \eqref{Lrho=0} for $\epsilon_1=0$ or $\epsilon_2=0$ follows by \eqref{Lorhoo=0}, \eqref{Lkrhoke0=0}, \eqref{Lkrhok0e=0}, and by a straightforward computation.
 \qed
 
 \vspace{\baselineskip}
 
 As done in Section \ref{M1M2}, we plan to apply (a corollary of) the implicit function theorem to prove the real analyticity of $(\theta^o[\cdot,\cdot], \theta_1[\cdot,\cdot], \theta_2[\cdot,\cdot])$. In order to do so, in the following technical Lemma \ref{dL}, we  study the regularity of $L$ and its partial differential.

\begin{lem}\label{dL}
The following statements hold.
\begin{enumerate}
\item[(i)] The map $L$ is real analytic from $]-\delta_1,\delta_1[\times]-\delta_2,\delta_2[\times C^{1,\alpha}(\partial\Omega^o)\times C^{1,\alpha}(\partial\Omega_1)_0\times C^{1,\alpha}(\partial\Omega_2)_0$ to $C^{1,\alpha}(\partial\Omega^o)\times C^{1,\alpha}(\partial\Omega_1)\times C^{1,\alpha}(\partial\Omega_2)$.
\item[(ii)] Let $(\bar\epsilon_1,\bar\epsilon_2,\bar\theta^o,\bar\theta_1,\bar\theta_2)\in ]-\delta_1,\delta_1[\times]-\delta_2,\delta_2[\times C^{1,\alpha}(\partial\Omega^o)\times C^{1,\alpha}(\partial\Omega_1)_0\times C^{1,\alpha}(\partial\Omega_2)_0$, then 
\begin{equation}\label{dL.eq1}
\partial_{(\theta^o,\theta_1,\theta_2)}L[\bar\epsilon_1,\bar\epsilon_2,\bar\theta^o,\bar\theta_1,\bar\theta_2]
\end{equation}
(the partial differential of $L$ with respect to the variable $(\theta^o,\theta_1,\theta_2)$ evaluated at $(\bar\epsilon_1,\bar\epsilon_2,\bar\theta^o,\bar\theta_1,\bar\theta_2)$) is an isomorphism from $C^{1,\alpha}(\partial\Omega^o)\times C^{1,\alpha}(\partial\Omega_1)_0\times C^{1,\alpha}(\partial\Omega_2)_0$ onto the subspace of  $C^{1,\alpha}(\partial\Omega^o)\times C^{1,\alpha}(\partial\Omega_1)\times C^{1,\alpha}(\partial\Omega_2)$ consisting of those triples $(\psi^o,\psi_1,\psi_2)$ such that
\begin{equation}\label{dL.eq2}
\int_{\partial\Omega^o}\psi^o\rho^o_h[\bar\epsilon_1,\bar\epsilon_2]\,d\sigma-\sum_{k=1}^2\int_{\partial\Omega_k}\psi_k\,\rho_{h,k}[\bar\epsilon_1,\bar\epsilon_2]\,d\sigma=0\quad\forall h\in\{1,2\}\,.
\end{equation}
\end{enumerate}
\end{lem}
\proof  
Statement (i) follows by the standard properties of integral  operators with real analytic kernels and with no singularity (see, e.g., Lanza de Cristoforis and the second author \cite{LaMu13}) and by classical mapping properties of layer potentials (cf., {\it e.g.}, Miranda \cite{Mi65}).

To prove statement (ii) we observe that  the partial differential \eqref{dL.eq1} is delivered by  
\[
\begin{split}
&\partial_{(\theta^o,\theta_1,\theta_2)}L^o[\bar\epsilon_1,\bar\epsilon_2,\bar\theta^o,\bar\theta_1,\bar\theta_2](\theta^o,\theta_1,\theta_2)=L^o[\bar\epsilon_1,\bar\epsilon_2,\theta^o,\theta_1,\theta_2]+f^o\\
&\partial_{(\theta^o,\theta_1,\theta_2)}L_k[\bar\epsilon_1,\bar\epsilon_2,\bar\theta^o,\bar\theta_1,\bar\theta_2](\theta^o,\theta_1,\theta_2)=L_k[\bar\epsilon_1,\bar\epsilon_2,\theta^o,\theta_1,\theta_2]\\
&\qquad -f_h+\int_{\partial\Omega^o}f^o\,\rho^o_h[\bar\epsilon_1,\bar\epsilon_2]\,d\sigma+\sum_{i=1}^2\int_{\partial\Omega_i}f_i\,\rho_{h,i}[\bar\epsilon_1,\bar\epsilon_2]\,d\sigma\qquad\forall k\in\{1,2\}\,,
\end{split}
\] 
for all $(\theta^o,\theta_1,\theta_2)\in  C^{1,\alpha}(\partial\Omega^o)\times C^{1,\alpha}(\partial\Omega_1)_0\times C^{1,\alpha}(\partial\Omega_2)_0$.  Then we take a triple $(\psi^o,\psi_1,\psi_2)$ in $C^{1,\alpha}(\partial\Omega^o)\times C^{1,\alpha}(\partial\Omega_1)\times C^{1,\alpha}(\partial\Omega_2)$ which satisfies condition \eqref{dL.eq2} and, by arguing as in the proof of Propositions \ref{Lee=0}, \ref{Le0=0}, \ref{L0e=0}, we verify that there exist  a unique triple $(\theta^o,\theta_1,\theta_2)\in  C^{1,\alpha}(\partial\Omega^o)\times C^{1,\alpha}(\partial\Omega_1)_0\times C^{1,\alpha}(\partial\Omega_2)_0$ such that 
\[
\partial_{(\theta^o,\theta_1,\theta_2)}L[\bar\epsilon_1,\bar\epsilon_2,\bar\theta^o,\bar\theta_1,\bar\theta_2](\theta^o,\theta_1,\theta_2)=(\psi^o,\psi_1,\psi_2)\,.
\]
Now the validity of the statement (ii) follows by the open mapping theorem and by Proposition \ref{Lperp}. 
\qed

\vspace{\baselineskip}

We now introduce in the following Lemma \ref{implicit} a technical corollary of the implicit function theorem for real analytic maps. For a proof we refer to Lanza de Cristoforis~\cite[Thm.~13]{La07}. 

\begin{lem}\label{implicit}
 Let $\mathcal{X}$, $\mathcal{Y}$, $\mathcal{Z}$, $\mathcal{Z}_{1}$ be Banach spaces. Let $\mathcal{O}$ be an open subset of $\mathcal{X}\times\mathcal{Y}$ such that $(\bar x,\bar y)\in\mathcal{O}$. Let $F$ be a real analytic map from $\mathcal{O}$ to $\mathcal{Z}$ such that $F(\bar x,\bar y)=0$. Let the partial differential $\partial_{y}F(\bar x,\bar y)$ with respect to the variable $y$ be an homeomorphism from $\mathcal{Y}$ onto its image $V\equiv \mathrm{Ran}(\partial_{y}F(\bar x,\bar y))$. Assume that there exists a closed subspace $V_{1}$ of $\mathcal{Z}$ such that $\mathcal{Z}=V\oplus V_{1}$. Let $\mathcal{O}_{1}$ be an open subset of $\mathcal{X}\times\mathcal{Y}\times\mathcal{Z}$ containing $(\bar x,\bar y,0)$ and such that $(x,y,F(x,y))$ and $(x,y,0)$ belong to $\mathcal{O}_{1}$ for all $(x,y)\in\mathcal{O}$. Let $G$ be a real analytic map from $\mathcal{O}_{1}$ to $\mathcal{Z}_{1}$ such that $G(x,y,F(x,y))=0$ for all $(x,y)\in\mathcal{O}$,    $G(x,y,0)=0$ for all $(x,y)\in\mathcal{O}$, and such that the partial differential $\partial_{z}G(\bar x,\bar y,0)$ is surjective onto $\mathcal{Z}_{1}$ and has kernel equal to $V$. Then there  exist an open neighbourhood $\mathcal{U}$ of $\bar x$ in $\mathcal{X}$, an open neighbourhood $\mathcal{V}$ of $\bar y$ in $\mathcal{Y}$ with $\mathcal{U}\times\mathcal{V}\subseteq\mathcal{O}$,  and a real analytic map $T$ from $\mathcal{U}$ to $\mathcal{V}$ such that the set of zeros of $F$ in $\mathcal{U}\times\mathcal{V}$ coincides with the graph of $T$.
\end{lem}

We are finally in the position to apply Lemma \ref{implicit} to equation $L[\epsilon_1,\epsilon_2,\theta^o,\theta_1,\theta_2]=0$ and prove that the triple $(\theta^o[\epsilon_1,\epsilon_2], \theta_1[\epsilon_1,\epsilon_2], \theta_2[\epsilon_1,\epsilon_2])$ depends analytically on  $(\epsilon_1,\epsilon_2)$.

\begin{prop}\label{thetaee}
The  function from $]-\delta_1,\delta_1[\times]-\delta_2,\delta_2[$ to  $C^{1,\alpha}(\partial\Omega^o)\times C^{1,\alpha}(\partial\Omega_1)_0\times C^{1,\alpha}(\partial\Omega_2)_0$ which takes $(\epsilon_1,\epsilon_2)$ to $(\theta^o[\epsilon_1,\epsilon_2], \theta_1[\epsilon_1,\epsilon_2], \theta_2[\epsilon_1,\epsilon_2])$ is real analytic. Moreover, the set of zeros of $L$ in $]-\delta_1,\delta_1[\times]-\delta_2,\delta_2[\times C^{1,\alpha}(\partial\Omega^o)\times C^{1,\alpha}(\partial\Omega_1)_0\times C^{1,\alpha}(\partial\Omega_2)_0$ coincides with the graph of $(\theta^o[\cdot,\cdot], \theta_1[\cdot,\cdot], \theta_2[\cdot,\cdot])$.
\end{prop}
\proof
Let $(\bar\epsilon_{1},\bar\epsilon_{2},\bar\theta^o,\bar\theta_{1},\bar\theta_{2})\in  ]-\delta_1,\delta_1[\times]-\delta_2,\delta_2[\times C^{1,\alpha}(\partial\Omega^o)\times C^{1,\alpha}(\partial\Omega_1)_0\times C^{1,\alpha}(\partial\Omega_2)_0$ be such that $L[\bar\epsilon_{1},\bar\epsilon_{2},\bar\theta^o,\bar\theta_{1},\bar\theta_{2}]=0$. Let $\mathcal{X}\equiv\mathbb{R}^2$, $\mathcal{Y}\equiv C^{1,\alpha}(\partial\Omega^o)\times C^{1,\alpha}(\partial\Omega_1)_0\times C^{1,\alpha}(\partial\Omega_2)_0$, $\mathcal{Z}\equiv C^{1,\alpha}(\partial\Omega^o)\times C^{1,\alpha}(\partial\Omega_1)\times C^{1,\alpha}(\partial\Omega_2)$, $\mathcal{Z}_1\equiv\mathbb{R}^2$, $\mathcal{O}\equiv  ]-\delta_1,\delta_1[\times]-\delta_2,\delta_2[\times C^{1,\alpha}(\partial\Omega^o)\times C^{1,\alpha}(\partial\Omega_1)_0\times C^{1,\alpha}(\partial\Omega_2)_0$. Let $F\equiv L$. Let $\bar x\equiv(\bar\epsilon_1,\bar \epsilon_2)$ and $\bar y\equiv (\bar\theta^o,\bar\theta_{1},\bar\theta_{2})$. Let $V$ be  the subspace of  $C^{1,\alpha}(\partial\Omega^o)\times C^{1,\alpha}(\partial\Omega_1)\times C^{1,\alpha}(\partial\Omega_2)$ consisting of the triples $(\psi^o,\psi_1,\psi_2)$  which satisfy the condition in \eqref{dL.eq2} with $\epsilon_1=\bar\epsilon_1$ and $\epsilon_2=\bar\epsilon_2$, let $V_1$ be the $2$-dimensional subspace of  $C^{1,\alpha}(\partial\Omega^o)\times C^{1,\alpha}(\partial\Omega_1)\times C^{1,\alpha}(\partial\Omega_2)$  generated by $(\rho^o_1[\bar\epsilon_1,\bar\epsilon_2],\rho_{1,1}[\bar\epsilon_1,\bar\epsilon_2],\rho_{1,2}[\bar\epsilon_1,\bar\epsilon_2] )$ and $(\rho^o_2[\bar\epsilon_1,\bar\epsilon_2],\rho_{2,1}[\bar\epsilon_1,\bar\epsilon_2],\rho_{2,2}[\bar\epsilon_1,\bar\epsilon_2] )$. Let 
$\mathcal{O}_1\equiv ]-\delta_1,\delta_1[\times]-\delta_2,\delta_2[\times C^{1,\alpha}(\partial\Omega^o)\times C^{1,\alpha}(\partial\Omega_1)_0\times C^{1,\alpha}(\partial\Omega_2)_0\times  C^{1,\alpha}(\partial\Omega^o)\times C^{1,\alpha}(\partial\Omega_1)\times C^{1,\alpha}(\partial\Omega_2)$. 
Let $G\equiv (G_1,G_2)$ be defined by
\[
G_h(\epsilon_1,\epsilon_2,\theta^o,\theta_1,\theta_2,\psi^o,\psi_1,\psi_2)
 \equiv\int_{\partial\Omega^o}\psi^o\rho^o_h[\epsilon_1,\epsilon_2]\,d\sigma-\sum_{k=1}^2\int_{\partial\Omega_k}\psi_k\,\rho_{h,k}[\epsilon_1,\epsilon_2]\,d\sigma 
\] 
for all $h\in\{1,2\}$  and for all $(\epsilon_1,\epsilon_2,\theta^o,\theta_1,\theta_2,\psi^o,\psi_1,\psi_2)\in\mathcal{O}_1$. Then Lemma \ref{implicit} implies that there exist an open neighbourhood of $\mathcal{U}$ of $(\bar\epsilon_1,\bar \epsilon_2)$ in $]-\delta_1,\delta_1[\times]-\delta_2,\delta_2[$, an open neighbourhood $\mathcal{V}$ of $ (\bar\theta^o,\bar\theta_{1},\bar\theta_{2})$  in $C^{1,\alpha}(\partial\Omega^o)\times C^{1,\alpha}(\partial\Omega_1)_0\times C^{1,\alpha}(\partial\Omega_2)_0$, and a real analytic map $T\equiv(T^o,T_1,T_2)$ from $\mathcal{U}$ to $\mathcal{V}$ such that the set of zeros of $L$ in $\mathcal{U}\times\mathcal{V}$ coincides with the graph of $T$. Then Propositions \ref{Lee=0}, \ref{Le0=0}, and \ref{L0e=0} imply that  $T[\epsilon_1,\epsilon_2]=(\theta^o[\epsilon_1,\epsilon_2],\theta_1[\epsilon_1,\epsilon_2], \theta_2[\epsilon_1,\epsilon_2])$ for all $(\epsilon_1,\epsilon_2)\in\mathcal{U}$ and the validity of the proposition follows.
\qed

\section{The auxiliary functions $H^{\Omega^o}_x$, $H^x_{\Omega_1}$, $H^x_{\Omega_2}$, and $H^x_{\tilde\Omega(\epsilon_2)}$}\label{H}
 
In the next Section \ref{rep}, we will exploit the results of Sections \ref{M1M2} and \ref{L} and the representation formula of Proposition \ref{ugO} to describe the dependence of the solution  $u_{\epsilon_1,\epsilon_2}$ of \eqref{dir} in terms of analytic functions of $\epsilon_1$, $\epsilon_2$ and of elementary functions of $\log|\epsilon_1|$ and $\log|\epsilon_1\epsilon_2|$. Before doing so, we introduce in this section the auxiliary functions $H^{\Omega^o}_x$, $H^x_{\Omega_1}$, $H^x_{\Omega_2}$, and $H^x_{\tilde\Omega(\epsilon_2)}$, which will play an important role in the description of the limit behaviour of $u_{\epsilon_1,\epsilon_2}$.  We note that $H^{\Omega^o}_x(y)$ is the difference between  the Dirichlet Green function in $\Omega^o$ and the fundamental solution $S(x-y)$ (see \eqref{Green}). Analogous relations hold for  $H^x_{\Omega_1}(y)$, $H^x_{\Omega_2}(y)$, and $H^x_{\tilde\Omega(\epsilon_2)}(y)$ in the exterior domains $\mathbb{R}^2\setminus\Omega_1$, $\mathbb{R}^2\setminus\Omega_2$, and $\mathbb{R}^2\setminus\tilde\Omega(\epsilon_2)$, respectively.

 \begin{prop}\label{HOox}
Let $x\in\Omega^o$ be fixed. Let $H^{\Omega^o}_x\in C^{1,\alpha}(\mathrm{cl}\Omega^o)$ be the solution of 
\[
\left\{
\begin{array}{ll}
\Delta H^{\Omega^o}_x=0&\text{in }\Omega^o\,,\\
H^{\Omega^o}_x(y)=S(x-y)&\forall y\in\partial\Omega^o\,.
\end{array}
\right.
\]
Then $v_{\Omega^o}[\rho^o_j[\epsilon_1,0]](x)=-H^{\Omega^o}_x(\epsilon_1p^j)$ and $v_{\Omega^o}[\rho^o_j[0,\epsilon_2]](x)=-H^{\Omega^o}_x(0)$ for all $(\epsilon_1,\epsilon_2)\in{]-\delta_1,\delta_1[}\times{]-\delta_2,\delta_2[}$ and for all $j\in\{1,2\}$.
 \end{prop}
 \proof
 Let $u\in C^{1,\alpha}(\mathrm{cl}\Omega^o)$ and $\Delta u=0$ in $\Omega^o$. Then by classical potential theory there exists $\mu\in C^{1,\alpha}(\partial\Omega^o)$ such that $u=w^+_{\Omega^o}[\mu]$ (cf.~Section \ref{prel}). Then, by the jump properties of the double layer potential (see \eqref{jump}), by standard properties of adjoint operators, and by Proposition \ref{Me0=0}, we have
\begin{equation}\label{HOox.eq1}
\begin{split}
&\int_{\partial\Omega^o} u_{|\partial\Omega^o}\rho^o_j[\epsilon_1,0]\,d\sigma=\int_{\partial\Omega^o} w^+_{\Omega^o}[\mu]_{|\partial\Omega^o}\rho^o_j[\epsilon_1,0]\,d\sigma\\
&\qquad=\int_{\partial\Omega^o} \Bigl[(\frac{1}{2}I_{\Omega^o}+W_{\Omega^o})\mu\Bigr]\rho^o_j[\epsilon_1,0]\,d\sigma =\int_{\partial\Omega^o} \mu\Bigl[(\frac{1}{2}I_{\Omega^o}+W^*_{\Omega^o})\rho^o_j[\epsilon_1,0]\Bigr]\,d\sigma\\
&\qquad=-\int_{\partial\Omega^o} \mu(y)\,\nu_{\Omega^o}(y)\cdot\nabla S(y-\epsilon_1p^j)\,d\sigma_y=-w^+_{\Omega^o}[\mu](\epsilon_1p^j)=-u(\epsilon_1p^j)\,.
\end{split}
\end{equation}
It follows that  
\[
v_{\Omega^o}[\rho^o_j[\epsilon_1,0]](x)=\int_{\partial\Omega^o}S(x-y)\rho^o_j[\epsilon_1,0](y)\,d\sigma=\int_{\partial\Omega^o}H^{\Omega^o}_{x|\partial\Omega^o}\,\rho^o_j[\epsilon_1,0]\,d\sigma=-H^{\Omega^o}_x(\epsilon_1p^j)\,.
\]

The proof of $v_{\Omega^o}[\rho^o_j[0,\epsilon_2]](x)=-H^{\Omega^o}_x(0)$ is similar. Indeed, for $u$ and $\mu$ as above we have 
\begin{equation}\label{HOox.eq2}
\begin{split}
&\int_{\partial\Omega^o} u_{|\partial\Omega^o}\rho^o_j[0,\epsilon_2]\,d\sigma=\int_{\partial\Omega^o} w^+_{\Omega^o}[\mu]_{|\partial\Omega^o}\rho^o_j[0,\epsilon_2]\,d\sigma\\
&\qquad=\int_{\partial\Omega^o} \Bigl[(\frac{1}{2}I_{\Omega^o}+W_{\Omega^o})\mu\Bigr]\rho^o_j[0,\epsilon_2]\,d\sigma =\int_{\partial\Omega^o} \mu\Bigl[(\frac{1}{2}I_{\Omega^o}+W^*_{\Omega^o})\rho^o_j[0,\epsilon_2]\Bigr]\,d\sigma\\
&\qquad=-\int_{\partial\Omega^o} \mu(y)\,\nu_{\Omega^o}(y)\cdot\nabla S(y)\,d\sigma_y=-w^+_{\Omega^o}[\mu](0)=-u(0)
\end{split}
\end{equation} (see also Proposition \ref{M0e=0}) and thus
\[
v_{\Omega^o}[\rho^o_j[0,\epsilon_2]](x)=\int_{\partial\Omega^o}S(x-y)\rho^o_j[0,\epsilon_2](y)\,d\sigma=\int_{\partial\Omega^o}H^{\Omega^o}_{x|\partial\Omega^o}\,\rho^o_j[0,\epsilon_2]\,d\sigma=-H^{\Omega^o}_x(0)\,.
\]
\qed

\vspace{\baselineskip}

\begin{prop}\label{HxOh}
Let $h\in\{1,2\}$ and $x\in\mathbb{R}^2\setminus\partial\Omega_h$ be fixed. Let $H^x_{\Omega_h}\in C^{1,\alpha}_{\mathrm{loc}}(\mathbb{R}^2\setminus\Omega_h)$ be the solution of 
\[
\left\{
\begin{array}{ll}
\Delta H^x_{\Omega_h}=0&\text{in }\mathbb{R}^2\setminus\mathrm{cl}\Omega_h\,,\\
H^x_{\Omega_h}(y)=S(x-y)&\forall y\in\partial\Omega_h\,,\\
\sup_{y\in\mathbb{R}^2\setminus\Omega_h}|H^x_{\Omega_h}(y)|<+\infty\,.
\end{array}
\right.
\]
Then 
\begin{equation}\label{HxOh.eq0}
v_{\Omega_h}[\rho_{j,h}[\epsilon_1,0]](x)=\delta_{j,h}\lim_{y\to\infty}H^x_{\Omega_h}(y)\qquad\forall j\in\{1,2\}\,.
\end{equation} 
If in addition $x\in\mathrm{cl}\Omega_h$, then we have 
\begin{equation}\label{HxOh.eq01}
v_{\Omega_h}[\rho_{j,h}[\epsilon_1,0]](x)=\delta_{j,h}\lim_{y\to\infty}H^0_{\Omega_h}(y)\qquad\forall j\in\{1,2\}\,.
\end{equation}
In particular,  $\lim_{y\to\infty}H^x_{\Omega_h}(y)=\lim_{y\to\infty}H^0_{\Omega_h}(y)$ for all $x\in\mathrm{cl}\Omega_h$ and all $h\in\{1,2\}$.
\end{prop}

\proof
We first prove \eqref{HxOh.eq0}. Let $u\in C^{1,\alpha}_\mathrm{loc}(\mathbb{R}^2\setminus\Omega_h)$, $\Delta u=0$ in $\mathbb{R}^2\setminus\mathrm{cl}\Omega_h$, and $\sup_{y\in\mathbb{R}^2\setminus\Omega_h}|u(y)|<+\infty$. Then, by classical potential theory there exists $\mu\in C^{1,\alpha}(\partial\Omega_h)$ such that $u=w^-_{\Omega_h}[\mu]+\lim_{y\to\infty}u(y)$ (cf. Folland \cite[Ch.~3]{Fo95}, see also Section \ref{prel}). Then by the jump properties of the double layer potential \eqref{jump}, by standard properties of adjoint operators, and by Proposition \ref{Me0=0}, we have
\begin{equation}\label{HxOh.eq1}
\begin{split}
\int_{\partial\Omega_h} u_{|\partial\Omega_h}\rho_{j,h}[\epsilon_1,0]&d\sigma=\int_{\partial\Omega_h} w^-_{\Omega_h}[\mu]_{|\partial\Omega_h}\rho_{j,h}[\epsilon_1,0]\,d\sigma+\delta_{j,h}\lim_{y\to\infty}u(y)\\
&=\int_{\partial\Omega_h} \Bigl[(-\frac{1}{2}I_{\Omega_h}+W_{\Omega_h})\mu\Bigr]\rho_{j,h}[\epsilon_1,0]\,d\sigma+\delta_{j,h}\lim_{y\to\infty}u(y)\\ &=\int_{\partial\Omega_h} \mu\Bigl[(-\frac{1}{2}I_{\Omega_h}+W^*_{\Omega_h})\rho_{j,h}[\epsilon_1,0]\Bigr]\,d\sigma+\delta_{j,h}\lim_{y\to\infty}u(y)=\delta_{j,h}\lim_{y\to\infty}u(y)\,.
\end{split}
\end{equation}
Thus 
\[
\begin{split}
&v_{\Omega_h}[\rho_{j,h}[\epsilon_1,0]](x)\\
&\quad=\int_{\partial\Omega_h}S(x-y)\rho_{j,h}[\epsilon_1,0](y)\,d\sigma_y=\int_{\partial\Omega_h}H^x_{\Omega_h}(y)\,\rho_{j,h}[\epsilon_1,0](y)\,d\sigma_y=\delta_{j,h}\lim_{y\to\infty}H^x_{\Omega_h}(y)\,.
\end{split}
\] 
To prove  \eqref{HxOh.eq01} we observe that, by  Proposition  \ref{Me0=0} and by the jump properties of the normal derivative of the single layer potential (cf.~\eqref{jump}), we have $\nu_{\Omega_h}\cdot\nabla v^+_{\Omega_h}[\rho_{j,h}[\epsilon_1,0]]_{|\partial\Omega_h}=0$. We deduce that $v_{\Omega_h}[\rho_{j,h}[\epsilon_1,0]]$ is constant on $\mathrm{cl}\Omega_h$ and the validity of statement (ii) follows.
\qed

\vspace{\baselineskip}

In the proof of Proposition \ref{HxOe} here below we exploit the following result of potential theory.

\begin{lem}\label{wcompleted}
Let $\epsilon_2\in]-\delta_2,\delta_2[\setminus\{0\}$ and let $h,k\in\{1,2\}$ with $h\neq k$. Then the operator from $C^{1,\alpha}(\partial\tilde\Omega(\epsilon_2))$ to itself which takes $\mu$ to the function defined by 
\begin{equation}\label{wcompleted.eq1}
\begin{split}
&\Bigl[(-\frac{1}{2}I_{\tilde\Omega(\epsilon_2)}+W_{\tilde\Omega(\epsilon_2)})\mu\Bigr](x)
\\
&\quad
 +\int_{\partial\Omega_h(1,\epsilon_2)}\mu\,d\sigma+(S(x-p^h)-S(x-p^k))\int_{\partial\Omega_k(1,\epsilon_2)}\mu\,d\sigma\qquad\forall x\in\partial\tilde\Omega(\epsilon_2)
\end{split}
\end{equation}
is a linear isomorphism. 
\end{lem}
\proof
By Proposition \ref{fredC} and by standard properties of Fredholm operators one verifies that the operator from $C^{1,\alpha}(\partial\tilde\Omega(\epsilon_2))$ to itself which takes a function $\mu$ to the function defined by  \eqref{wcompleted.eq1} is Fredholm of index $0$. Thus, in order to show that it is an isomorphism it suffices to show that $\mu=0$ when 
 \begin{equation}\label{wcompleted.eq2}
\begin{split}
&\Bigl[(-\frac{1}{2}I_{\tilde\Omega(\epsilon_2)}+W_{\tilde\Omega(\epsilon_2)})\mu\Bigr](x)
\\
&\quad
 +\int_{\partial\Omega_h(1,\epsilon_2)}\mu\,d\sigma+(S(x-p^h)-S(x-p^k))\int_{\partial\Omega_k(1,\epsilon_2)}\mu\,d\sigma=0\qquad\forall x\in\partial\tilde\Omega(\epsilon_2)\,.
\end{split}
\end{equation}
If $\mu$ satisfies equation \eqref{wcompleted.eq2}, then by the jump properties of the double layer potential (cf.~\eqref{jump}) we have 
\begin{equation}\label{wcompleted.eq3}
w^-_{\tilde\Omega(\epsilon_2)}[\mu](x)=-\int_{\partial\Omega_h(1,\epsilon_2)}\mu\,d\sigma-(S(x-p^h)-S(x-p^k))\int_{\partial\Omega_k(1,\epsilon_2)}\mu\,d\sigma
\end{equation}
for all $x\in\partial\tilde\Omega(\epsilon_2)$. We observe that both the left and the right hand side of \eqref{wcompleted.eq3} define functions which are bounded in $\mathbb{R}^2\setminus\tilde\Omega(\epsilon_2)$. Accordingly, the uniqueness properties of the solution of the exterior Dirichlet problem (cf., {\it e.g.}, Folland \cite[Chap.~2]{Fo95}) implies that equality \eqref{wcompleted.eq3} holds for all $x\in\mathbb{R}^2\setminus\tilde\Omega(\epsilon_2)$. Then, by the decay properties of $w^-_{\tilde\Omega(\epsilon_2)}[\mu](x)$ and of $S(x-p^h)-S(x-p^k)$ as $x\to\infty$ we deduce that 
\begin{equation}\label{wcompleted.eq4}
\int_{\partial\Omega_h(1,\epsilon_2)}\mu\,d\sigma=0\,.
\end{equation}
Now we observe that by equality \eqref{nojump} and by the divergence theorem we have
\begin{equation}\label{wcompleted.eq5}
\begin{split}
&\int_{\partial\Omega_h(1,\epsilon_2)}\nu_{\Omega_h(1,\epsilon_2)}(x)\cdot\nabla w^-_{\tilde\Omega(\epsilon_2)}[\mu](x)\,d\sigma_x
\\
&\qquad\qquad=\int_{\partial\Omega_h(1,\epsilon_2)}\nu_{\Omega_h(1,\epsilon_2)}(x)\cdot\nabla w^+_{\tilde\Omega(\epsilon_2)}[\mu](x)\,d\sigma_x=0\,.
\end{split}
\end{equation}
Moreover, by the definition of the double layer potential and by equalities $w_{\Omega_h(1,\epsilon_2)}[1](p^k)=0$ and $w_{\Omega_h(1,\epsilon_2)}[1](p^h)=1$ (cf.~Section \ref{prel}, see also Folland \cite[Chap.~3]{Fo95}) we have 
\begin{equation}\label{wcompleted.eq6}
\begin{split}
&\int_{\partial\Omega_h(1,\epsilon_2)}\nu_{\Omega_h(1,\epsilon_2)}(x)\cdot\nabla (S(x-p^h)-S(x-p^k))\,d\sigma_x=w_{\Omega_h(1,\epsilon_2)}[1](p^h)-w_{\Omega_h(1,\epsilon_2)}[1](p^k)=1\,.
\end{split}
\end{equation}
Hence, by equalities \eqref{wcompleted.eq3}, \eqref{wcompleted.eq5}, and \eqref{wcompleted.eq6} we deduce that 
\begin{equation}\label{wcompleted.eq7}
\int_{\partial\Omega_k(1,\epsilon_2)}\mu\,d\sigma=0\,.
\end{equation}
Then, by equalities \eqref{wcompleted.eq3}, \eqref{wcompleted.eq4}, and \eqref{wcompleted.eq7}, and by Lemma \ref{ker} it follows that $\mu=0$. Our proof is now completed.
\qed

\vspace{\baselineskip}

We observe here that Lemma \ref{wcompleted} implies that
\begin{equation}\label{nuu1=-nuu2}
\int_{\partial\Omega_1(1,\epsilon_2)}\nu_{\Omega_1(1,\epsilon_2)}(y)\cdot\nabla u(y)\,d\sigma_y=-\int_{\partial\Omega_2(1,\epsilon_2)}\nu_{\Omega_2(1,\epsilon_2)}(y)\cdot\nabla u(y)\,d\sigma_y
\end{equation}
for all $u\in C^{1,\alpha}_\mathrm{loc}(\mathbb{R}^2\setminus\tilde\Omega(\epsilon_2))$ such that $\Delta u=0$ in $\mathbb{R}^2\setminus\mathrm{cl}\tilde\Omega(\epsilon_2)$ and $\sup_{y\in\mathbb{R}^2\setminus\tilde\Omega(\epsilon_2)}|u(y)|<+\infty$ (see also Folland \cite[Chap.~2]{Fo95}). 

\begin{prop}\label{HxOe}
Let $\epsilon_2\in]-\delta_2,\delta_2[\setminus\{0\}$ be fixed. For each $x\in\mathbb{R}^2\setminus\partial\tilde\Omega(\epsilon_2)$ let $H^x_{\tilde\Omega(\epsilon_2)}\in C^{1,\alpha}_{\mathrm{loc}}(\mathbb{R}^2\setminus\tilde\Omega(\epsilon_2))$ be the solution of  
\begin{equation}\label{HxOe.eq0}
\left\{
\begin{array}{ll}
\Delta H^x_{\tilde\Omega(\epsilon_2)}=0&\text{in }\mathbb{R}^2\setminus\mathrm{cl}\tilde\Omega(\epsilon_2)\,,\\
H^x_{\tilde\Omega(\epsilon_2)}(y)=S(x-y)&\forall y\in\partial\tilde\Omega(\epsilon_2)\,,\\
\sup_{y\in\mathbb{R}^2\setminus\tilde\Omega(\epsilon_2)}|H^x_{\tilde\Omega(\epsilon_2)}(y)|<+\infty\,.
\end{array}
\right.
\end{equation}
Let $j\in\{1,2\}$. Let $\tilde\rho_j\in C^{0,\alpha}(\partial\tilde\Omega(\epsilon_2))$ be defined as in \eqref{tilderhoi}. Let $H^{j,i}_{\tilde\Omega(\epsilon_2)}\in\mathbb{R}$ be defined by 
\begin{equation}\label{HxOe.eq1}
H^{j,i}_{\tilde\Omega(\epsilon_2)}\equiv v_{\tilde\Omega(\epsilon_2)}[\tilde\rho_j](p^i)\qquad\forall i\in\{1,2\}\,.
\end{equation}
Let $h,k\in\{1,2\}$ with $h\neq k$. Then 
\begin{equation}\label{HxOe.eq2}
\begin{split}
&v_{\Omega_h}[\rho_{j,h}[0,\epsilon_2]](\xi)+\frac{\log|\epsilon_2|}{2\pi}\delta_{j,h}+\int_{\partial\Omega_k}S(p^h-p^k+\epsilon_2(\xi-\eta))\rho_{j,k}[0,\epsilon_2](\eta)\,d\sigma_\eta\\
&=\lim_{y\to\infty}H_{\tilde\Omega(\epsilon_2)}^{p^h+\epsilon_2\xi}(y)+\left(H^{j,k}_{\tilde\Omega(\epsilon_2)}-H^{j,h}_{\tilde\Omega(\epsilon_2)}\right)\int_{\partial\Omega_k(1,\epsilon_2)}\nu_{\Omega_k(1,\epsilon_2)}(y)\cdot\nabla_y H_{\tilde\Omega(\epsilon_2)}^{p^h+\epsilon_2\xi}(y)\,d\sigma_y\\
&=\lim_{y\to\infty}H_{\tilde\Omega(\epsilon_2)}^{p^h+\epsilon_2\xi}(y)+\left(H^{j,h}_{\tilde\Omega(\epsilon_2)}-H^{j,k}_{\tilde\Omega(\epsilon_2)}\right)\int_{\partial\Omega_h(1,\epsilon_2)}\nu_{\Omega_h(1,\epsilon_2)}(y)\cdot\nabla_y H_{\tilde\Omega(\epsilon_2)}^{p^h+\epsilon_2\xi}(y)\,d\sigma_y
\end{split}
\end{equation}
for all $\xi\in\mathbb{R}^2$ such that $p^h+\epsilon_2\xi\notin \partial\tilde\Omega(\epsilon_2)$.
In addition, if $\xi\in \mathrm{cl}\Omega_h$, then 
\begin{equation}\label{HxOe.eq1.1}
v_{\Omega_h}[\rho_{j,h}[0,\epsilon_2]](\xi)+\frac{\log|\epsilon_2|}{2\pi}\delta_{j,h}+\int_{\partial\Omega_k}S(p^h-p^k+\epsilon_2(\xi-\eta))\rho_{j,k}[0,\epsilon_2](\eta)\,d\sigma_\eta=H^{j,h}_{\tilde\Omega(\epsilon_2)}
\end{equation}
and if $\xi\in (p^k-p^h)/\epsilon_2+\mathrm{cl}\Omega_k$, then 
\begin{equation}\label{HxOe.eq1.2}
v_{\Omega_h}[\rho_{j,h}[0,\epsilon_2]](\xi)+\frac{\log|\epsilon_2|}{2\pi}\delta_{j,h}+\int_{\partial\Omega_k}S(p^h-p^k+\epsilon_2(\xi-\eta))\rho_{j,k}[0,\epsilon_2](\eta)\,d\sigma_\eta=H^{j,k}_{\tilde\Omega(\epsilon_2)}\,.
\end{equation}
\end{prop}
\proof
Let $u\in C^{1,\alpha}_\mathrm{loc}(\mathbb{R}^2\setminus\tilde\Omega(\epsilon_2))$, $\Delta u=0$ in $\mathbb{R}^2\setminus\mathrm{cl}\tilde\Omega(\epsilon_2)$, and $\sup_{y\in\mathbb{R}^2\setminus\tilde\Omega(\epsilon_2)}|u(y)|<+\infty$. Then, by classical potential theory there exists $\mu\in C^{1,\alpha}(\partial\tilde\Omega(\epsilon_2))$ such that 
\[
u(x)=w^-_{\tilde\Omega(\epsilon_2)}[\mu](x)+\int_{\partial\Omega_h(1,\epsilon_2)}\mu\,d\sigma+(S(x-p^h)-S(x-p^k))\int_{\partial\Omega_k(1,\epsilon_2)}\mu\,d\sigma\qquad\forall x\in\mathbb{R}^2\setminus\tilde\Omega(\epsilon_2)
\] 
(cf.~Lemma \ref{wcompleted}). Then, a computation based on the divergence theorem, on equality \eqref{nojump}, and on equality $w_{\Omega_k(1,\epsilon_2)}[1](p^k)=1$, shows that
\[
\begin{split}
\int_{\partial\Omega_k(1,\epsilon_2)}\nu_{\Omega_k(1,\epsilon_2)}\cdot\nabla u\, d\sigma& =-\int_{\partial\Omega_k(1,\epsilon_2)}\nu_{\Omega_k(1,\epsilon_2)}\cdot\nabla S(x-p^k)\, d\sigma\int_{\partial\Omega_k(1,\epsilon_2)}\mu\, d\sigma\\
&=-w_{\Omega_k(1,\epsilon_2)}[1](p^k)\int_{\partial\Omega_k(1,\epsilon_2)}\mu\, d\sigma=-\int_{\partial\Omega_k(1,\epsilon_2)}\mu\, d\sigma\,.
\end{split}
\]
Hence,  by the jump properties of the double layer potential \eqref{jump} and by the decay at $\infty$ of $w^-_{\Omega_h}[\mu](x)$ and $S(x-p^k)-S(x-p^h)$  we obtain that
\[
u(x)=w^-_{\tilde\Omega(\epsilon_2)}[\mu](x)+\lim_{y\to\infty}u(y)+(S(x-p^k)-S(x-p^h))\int_{\partial\Omega_k(1,\epsilon_2)}\nu_{\Omega_k(1,\epsilon_2)}\cdot\nabla u\,d\sigma\quad\forall x\in\mathbb{R}^2\setminus\tilde\Omega(\epsilon_2)\,.
\]
 Now let $\tilde\rho_j\in C^{0,\alpha}(\partial\tilde\Omega(\epsilon_2))$ be defined as in \eqref{tilderhoi}. Then by the jump properties of the double layer potential \eqref{jump}, by the definition of the single layer potential (cf.~Section \ref{prel}), by standard properties of adjoint operators, and by Proposition \ref{M0e=0}, we have
\begin{equation}\label{HxOe.eq3}
\begin{split}
&\int_{\partial\tilde\Omega(\epsilon_2)} u_{|\partial\tilde\Omega(\epsilon_2)}\,\tilde\rho_j\,d\sigma\\
&=\int_{\partial\tilde\Omega(\epsilon_2)} w^-_{\tilde\Omega(\epsilon_2)}[\mu]_{|\partial\tilde\Omega(\epsilon_2)}\,\tilde\rho_j\,d\sigma+\lim_{y\to\infty}u(y)\int_{\partial\tilde\Omega(\epsilon_2)} \tilde\rho_j\,d\sigma\\
&\quad +\left(\int_{\partial\tilde\Omega(\epsilon_2)}S(y-p^k)\,\tilde\rho_j(y)\,d\sigma_y-\int_{\partial\tilde\Omega(\epsilon_2)}S(y-p^h)\,\tilde\rho_j(y)\,d\sigma_y\right)\int_{\partial\Omega_k(1,\epsilon_2)}\nu_{\Omega_k(1,\epsilon_2)}\cdot\nabla u\,d\sigma\\
&=\int_{\partial\tilde\Omega(\epsilon_2)} \Bigl[(-\frac{1}{2}I_{\tilde\Omega(\epsilon_2)}+W_{\tilde\Omega(\epsilon_2)})\mu\Bigr]\tilde\rho_j\,d\sigma+\lim_{y\to\infty}u(y)\\
&\quad +\left(v_{\tilde\Omega(\epsilon_2)}[\tilde\rho_j](p^k)-v_{\tilde\Omega(\epsilon_2)}[\tilde\rho_j](p^h)\right)\int_{\partial\Omega_k(1,\epsilon_2)}\nu_{\Omega_k(1,\epsilon_2)}\cdot\nabla u\,d\sigma
\\ 
&=\int_{\partial\tilde\Omega(\epsilon_2)} \mu\Bigl[(-\frac{1}{2}I_{\tilde\Omega(\epsilon_2)}+W^*_{\tilde\Omega(\epsilon_2)})\tilde\rho_j\Bigr]\,d\sigma+\lim_{y\to\infty}u(y)\\
&\quad +\left(v_{\tilde\Omega(\epsilon_2)}[\tilde\rho_j](p^k)-v_{\tilde\Omega(\epsilon_2)}[\tilde\rho_j](p^h)\right)\int_{\partial\Omega_k(1,\epsilon_2)}\nu_{\Omega_k(1,\epsilon_2)}\cdot\nabla u\,d\sigma
\\
&=\lim_{y\to\infty}u(y)+\left(v_{\tilde\Omega(\epsilon_2)}[\tilde\rho_j](p^k)-v_{\tilde\Omega(\epsilon_2)}[\tilde\rho_j](p^h)\right)\int_{\partial\Omega_k(1,\epsilon_2)}\nu_{\Omega_k(1,\epsilon_2)}\cdot\nabla u\,d\sigma\,.
\end{split}
\end{equation}
Then, by the rule of change of variables in integrals and by \eqref{HxOe.eq0} we deduce that
\[
\begin{split}
&v_{\Omega_h}[\rho_{j,h}[0,\epsilon_2]](\xi)+\frac{\log|\epsilon_2|}{2\pi}\delta_{j,h}+\int_{\partial\Omega_k}S(p^h-p^k+\epsilon_2(\xi-\eta))\rho_{j,k}[0,\epsilon_2](\eta)\,d\sigma_\eta\\
&\quad =v_{\tilde\Omega(\epsilon_2)}[\tilde\rho^j](p^h+\epsilon_2\xi)\\
&\quad=\int_{\partial\tilde\Omega(\epsilon_2)}S(p^h+\epsilon_2\xi-y)\tilde\rho^j(y)\,d\sigma_y\\
&\quad =\lim_{y\to\infty}H^{p^h+\epsilon_2\xi}_{\tilde\Omega(\epsilon_2)}(y)\\
&\qquad +\left(v_{\tilde\Omega(\epsilon_2)}[\tilde\rho_j](p^k)-v_{\tilde\Omega(\epsilon_2)}[\tilde\rho_j](p^h)\right)\int_{\partial\Omega_k(1,\epsilon_2)}\nu_{\Omega_k(1,\epsilon_2)}(y)\cdot\nabla_y H^{p^h+\epsilon_2\xi}_{\tilde\Omega(\epsilon_2)}(y)\,d\sigma_y
\end{split}
\]
for all $\xi\in\mathbb{R}^2$ such that $p^h+\epsilon_2\xi\notin \partial\tilde\Omega(\epsilon_2)$.  It follows that the first equality in \eqref{HxOe.eq2} holds with $H^{j,k}_{\tilde\Omega(\epsilon_2)}$ and $H^{j,h}_{\tilde\Omega(\epsilon_2)}$ as in \eqref{HxOe.eq1}. Then, by \eqref{nuu1=-nuu2} one deduces the validity of the second equality in  \eqref{HxOe.eq2}. To prove \eqref{HxOe.eq1.1} and \eqref{HxOe.eq1.2} we observe that,  by  Proposition  \ref{Me0=0} and by the jump properties of the single layer potential \eqref{jump}, we have $\nu_{\tilde\Omega(\epsilon_2)}\cdot\nabla v^+_{\tilde\Omega(\epsilon_2)}[\tilde\rho_j]_{|\partial\tilde\Omega(\epsilon_2)}=0$. Thus  $v^+_{\tilde\Omega(\epsilon_2)}[\tilde\rho_j]$ is constant in $\mathrm{cl}\Omega_h(1,\epsilon_2)$ and in  $\mathrm{cl}\Omega_k(1,\epsilon_2)$ and the validity of \eqref{HxOe.eq1.1} and \eqref{HxOe.eq1.2}  follows by \eqref{HxOe.eq1} and by a straightforward computation based on the rule of change of variables in integrals.
\qed

\section{Representation of $u_{\epsilon_1,\epsilon_2}$ in terms of analytic maps}\label{rep}

In this section, we prove our main Theorem \ref{ue1e2} on the representation of $u_{\epsilon_1,\epsilon_2}$ in terms of real analytic maps and known functions. We will do so by exploiting the representation formula of Proposition \ref{ugO}, the real analyticity results of Propositions \ref{rhoee} and \ref{thetaee}, and the auxiliary functions of Sections \ref{H}.

In the following Propositions \ref{U}--\ref{Lambda} we introduce the functions $U[\epsilon_1,\epsilon_2]$ and $V[\epsilon_1,\epsilon_2]$, the vector $F[\epsilon_1,\epsilon_2]$, and the matrices $R[\epsilon_1,\epsilon_2]$ and $\Lambda(\epsilon_1,\epsilon_2)$ which we exploit to write $u_{\epsilon_1,\epsilon_2}$ and $u_{\epsilon_1,\epsilon_2}(\epsilon_1p^1 +\epsilon_1\epsilon_2\,\cdot\,)$ in terms of real analytic maps (cf.~Theorem \ref{ue1e2}).

\begin{prop}\label{U}
For each $(\epsilon_1,\epsilon_2)\in ]-\delta_1,\delta_1[\times]-\delta_2,\delta_2[$ there exists a unique function $U[\epsilon_1,\epsilon_2]$ in  $C^{1,\alpha}(\mathrm{cl}\Omega(\epsilon_1,\epsilon_2))$ such that 
\[
U[\epsilon_1,\epsilon_2](x)= w^+_{\Omega^o}[\theta^o[\epsilon_1,\epsilon_2]](x)+\epsilon_1\epsilon_2\sum_{k=1}^2\int_{\partial\Omega_k}\nu_{\Omega_k}(\eta)\cdot\nabla S(x-\epsilon_1p^k-\epsilon_1\epsilon_2\eta)\theta_{k}[\epsilon_1,\epsilon_2](\eta)\,d\sigma_\eta
\] 
for all $x\in\mathrm{cl}\Omega^o\setminus(\mathrm{cl}\Omega_1(\epsilon_1,\epsilon_2)\cup\mathrm{cl}\Omega_2(\epsilon_1,\epsilon_2))$. Moreover, the following statements hold.
\begin{enumerate}
\item[(i)] Let $\Omega_M$ be an open subset of $\Omega^o$ such that $0\notin\mathrm{cl}\Omega_M$. Let $\delta_M\in]0,\delta_1]$ be such that $\mathrm{cl}\Omega_M\cap\mathrm{cl}\Omega_k(\epsilon_1,\epsilon_2)=\emptyset$ for all $(\epsilon_1,\epsilon_2)\in]-\delta_M,\delta_M[\times]-\delta_2,\delta_2[$ and for all $k\in\{1,2\}$. Then there exists a real analytic map $U^M$ from $]-\delta_M,\delta_M[\times]-\delta_2,\delta_2[$ to $C^{1,\alpha}(\mathrm{cl}\Omega_M)$ such that 
\begin{equation}\label{U.eq1}
U[\epsilon_1,\epsilon_2](x)=u^o(x)+\epsilon_1\epsilon_2\,U^M[\epsilon_1,\epsilon_2](x)\quad\forall x\in\mathrm{cl}\Omega_M\,,\;(\epsilon_1,\epsilon_2)\in]-\delta_M,\delta_M[\times]-\delta_2,\delta_2[
\end{equation}
where $u^o\in C^{1,\alpha}(\mathrm{cl}\Omega^o)$ is the unique solution of 
\[
\left\{
\begin{array}{ll}
\Delta u^o=0&\text{in }\Omega^o\,,\\
u^o=f^o&\text{on }\partial\Omega^o\,.
\end{array}
\right.
\]
\item[(ii)]  Let $h,k\in\{1,2\}$ and $h\neq k$. Let $\Omega_m$ be an open bounded subset of $\mathbb{R}^2\setminus\mathrm{cl}\Omega_h$. Let $\delta_m\in]0,\delta_1]$ be such that $\epsilon_1p^h+\epsilon_1\epsilon_2\mathrm{cl}\Omega_m\subseteq\Omega^o$ and $(\epsilon_1p^h+\epsilon_1\epsilon_2\mathrm{cl}\Omega_m)\cap\mathrm{cl}\Omega_k(\epsilon_1,\epsilon_2)=\emptyset$ for all $(\epsilon_1,\epsilon_2)\in]-\delta_m,\delta_m[^2$. Then there exists a real analytic map $U^m_h$ from  $]-\delta_m,\delta_m[^2$ to $C^{1,\alpha}(\mathrm{cl}\Omega_m)$ such that 
\begin{equation}\label{U.eq2}
U[\epsilon_1,\epsilon_2](\epsilon_1p^h+\epsilon_1\epsilon_2\xi)=U^m_h[\epsilon_1,\epsilon_2](\xi)\qquad\forall \xi\in\mathrm{cl}\Omega_m\,,\; (\epsilon_1,\epsilon_2)\in(]-\delta_m,\delta_m[\setminus\{0\})^2\,.
\end{equation}
Moreover,
\begin{equation}\label{U.eq3}
U^m_h[\epsilon_1,0](\xi)=u^o(\epsilon_1 p^h)+u_h(\xi)-\lim_{\eta\to\infty} u_h(\eta)\quad\forall \xi\in\mathrm{cl}\Omega_m\,,\;\epsilon_1\in]-\delta_m,\delta_m[\,
\end{equation}
where $u_h\in C^{1,\alpha}_{\mathrm{loc}}(\mathbb{R}^2\setminus\Omega_h)$ is the solution of 
\[
\left\{
\begin{array}{ll}
\Delta u_h=0&\text{in }\mathbb{R}^2\setminus\mathrm{cl}\Omega_h\,,\\
u_{h}=f_h&\text{on }\partial\Omega_h\,,\\
\sup_{\eta\in\mathbb{R}^2\setminus\Omega_h}|u_h(\eta)|<+\infty\,,
\end{array}
\right.
\]
and 
\begin{equation}\label{U.eq4}
U_h^m[0,\epsilon_2](\xi)=u^o(0)+\tilde u(p^h+\epsilon_2\xi)+\tilde w(p^h+\epsilon_2\xi)-\lim_{y\to\infty}\tilde u(y)
\end{equation}
for all $\xi\in\mathrm{cl}\Omega_m$ and all $\epsilon_2\in]-\delta_m,\delta_m[\setminus\{0\}$,
 where $\tilde u\in C^{1,\alpha}_{\mathrm{loc}}(\mathbb{R}^2\setminus\tilde\Omega(\epsilon_2))$ is the solution of 
\begin{equation}\label{tildeu}
\left\{
\begin{array}{ll}
\Delta \tilde u=0&\text{in }\mathbb{R}^2\setminus\mathrm{cl}\tilde\Omega(\epsilon_2)\,,\\
\tilde u=\tilde f&\text{on }\partial\tilde\Omega(\epsilon_2)\,,\\
\sup_{y\in\mathbb{R}^2\setminus\tilde\Omega(\epsilon_2)}|\tilde u(y)|<+\infty\,,
\end{array}
\right.
\end{equation}
with  $\tilde{f}(x)\equiv
f_j((x-p^j)/\epsilon_2)$ for all $j\in\{1,2\}$ and $x\in\partial\Omega_j(1,\epsilon_2)$, and where $\tilde w\in C^{1,\alpha}_{\mathrm{loc}}(\mathbb{R}^2\setminus\tilde\Omega(\epsilon_2))$ is the solution of 
\[
\left\{
\begin{array}{ll}
\Delta \tilde w=0&\text{in }\mathbb{R}^2\setminus\mathrm{cl}\tilde\Omega(\epsilon_2)\,,\\
\tilde w(x)=(H^{i,j}_{\tilde\Omega(\epsilon_2)}-H^{i,i}_{\tilde\Omega(\epsilon_2)})\int_{\partial\Omega_j(1,\epsilon_2)}\nu_{\Omega_j(1,\epsilon_2)}\cdot\nabla\tilde u\,d\sigma&\forall i,j\in\{1,2\}\,,\; i\neq j\,, x\in\partial\Omega_i(1,\epsilon_2)\,,\\
\sup_{y\in\mathbb{R}^2\setminus\tilde\Omega(\epsilon_2)}|\tilde w(y)|<+\infty\,.
\end{array}
\right.
\]
\end{enumerate}
\end{prop}
\proof
We first consider statement (i). We observe that by Propositions \ref{Le0=0} and \ref{L0e=0} we have $\theta^o[0,\epsilon_2]=\theta^o[\epsilon_1,0]=\mu^o$ for all $(\epsilon_1,\epsilon_2)\in]-\delta_1,\delta_1[\times]-\delta_2,\delta_2[$, where $\mu^o\in C^{1,\alpha}(\partial\Omega^o)$ is the unique solution of $(\frac{1}{2}I_{\Omega^o}+W_{\Omega^o})\mu^o=f^o$. By standard properties of real analytic maps it follows that there is a real analytic map $\Theta^o$ from $]-\delta_1,\delta_1[\times]-\delta_2,\delta_2[$ to $C^{1,\alpha}(\partial\Omega^o)$ such that $\theta^o[\epsilon_1,\epsilon_2]=\mu^o+\epsilon_1\epsilon_2\,\Theta^o[\epsilon_1,\epsilon_2]$ for all $(\epsilon_1,\epsilon_2)\in]-\delta_1,\delta_1[\times]-\delta_2,\delta_2[$. Since $w^+_{\Omega^o}[\mu^o]=u^o$ by the jump formula \eqref{jump} and by the uniqueness of the solution of the Dirichlet problem, we deduce that 
\[
w^+_{\Omega^o}[\theta^o[\epsilon_1,\epsilon_2]]=u^o+\epsilon_1\epsilon_2\,w^+_{\Omega^o}[\Theta^o[\epsilon_1,\epsilon_2]]\,.
\]  
Then we define 
\[
U^M[\epsilon_1,\epsilon_2](x)\equiv w^+_{\Omega^o}[\Theta^o[\epsilon_1,\epsilon_2]](x)+\sum_{k=1}^2\int_{\partial\Omega_k}\nu_{\Omega_k}(\eta)\cdot\nabla S(x-\epsilon_1p^k-\epsilon_1\epsilon_2\eta)\theta_{k}[\epsilon_1,\epsilon_2](\eta)\,d\sigma_\eta
\] 
for all $x\in\mathrm{cl}\Omega_M$ and for all $(\epsilon_1,\epsilon_2)\in]-\delta_M,\delta_M[\times]-\delta_2,\delta_2[$. One readily verifies the validity of \eqref{U.eq1}. In addition,  by the standard properties of integral  operators with real analytic kernels and with no singularity (see, e.g., Lanza de Cristoforis and the second author \cite{LaMu13}), by the classical mapping properties of layer potentials (cf., {\it e.g.}, Miranda \cite{Mi65}), and by Proposition \ref{thetaee} one verifies that the map from $]-\delta_M,\delta_M[\times]-1,1[$ to $C^{1,\alpha}(\mathrm{cl}\Omega_M)$ which takes $(\epsilon_1,\epsilon_2)$ to  $U^M[\epsilon_1,\epsilon_2]$ is real analytic. 

We now prove statement (ii). We define 
\[
\begin{split}
U^m_h[\epsilon_1,\epsilon_2](\xi)&\equiv w_{\Omega^o}[\theta^o[\epsilon_1,\epsilon_2]](\epsilon_1p^h+\epsilon_1\epsilon_2\xi)-w^-_{\Omega_h}[\theta_h[\epsilon_1,\epsilon_2]](\xi)\\
&\quad +\epsilon_2\int_{\partial\Omega_k}\nu_{\Omega_k}(\eta)\cdot\nabla S(p^h-p^k+\epsilon_2(\xi-\eta))\theta_k[\epsilon_1,\epsilon_2](\eta)\,d\sigma_\eta\qquad\forall\xi\in\mathrm{cl}\Omega_m
\end{split}
\]
for all $(\epsilon_1,\epsilon_2)\in]-\delta_m,\delta_m[^2$. Then, by the standard properties of integral  operators with real analytic kernels and with no singularity (see, e.g., Lanza de Cristoforis and the second author \cite{LaMu13}) and by the classical mapping properties of layer potentials (cf., {\it e.g.}, Miranda \cite{Mi65}) we verify that the map which takes $(\epsilon_1,\epsilon_2)$ to $U^m_h[\epsilon_1,\epsilon_2]$ is real analytic from $]-\delta_m,\delta_m[^2$ to $C^{1,\alpha}(\mathrm{cl}\Omega_m)$. The validity of equality \eqref{U.eq2} can be deduced by a straightforward computation based on the rule of change of variables in integrals.  We now verify \eqref{U.eq3}. A straightforward computation shows that
\begin{equation}\label{U.eq5}
U^m_h[\epsilon_1,0](\xi)= w_{\Omega^o}[\theta^o[\epsilon_1,0]](\epsilon_1p^h)-w^-_{\Omega_h}[\theta_h[\epsilon_1,0]](\xi)\qquad\forall\xi\in\mathrm{cl}\Omega_m\,.
\end{equation}
Then we observe that by Proposition \ref{Le0=0} and by the jump formulae \eqref{jump} we have 
\begin{equation}\label{U.eq6}
w_{\Omega^o}[\theta^o[\epsilon_1,0]](\epsilon_1p^h)=u^o(\epsilon_1p_h)\,.
\end{equation}
In addition, by Proposition \ref{Le0=0} and by the jump formulae \eqref{jump}, we have 
\[
-w^-_{\Omega_h}[\theta_h[\epsilon_1,0]]_{|\partial\Omega_h}=f_h-\int_{\partial\Omega_h}f_h\,\rho_{h,h}[\epsilon_1,0]\,d\sigma\,.
\]
Accordingly, equality \eqref{HxOh.eq1} implies that
\[
-w^-_{\Omega_h}[\theta_h[\epsilon_1,0]]_{|\partial\Omega_h}=f_h-\lim_{y\to\infty}u_h(y)
\]
and by the uniqueness of the solution of the exterior Dirichlet problem we deduce that
\begin{equation}\label{U.eq7}
-w^-_{\Omega_h}[\theta_h[\epsilon_1,0]]=u_{h}-\lim_{y\to\infty}u_h(y)\qquad\text{on }\mathbb{R}^2\setminus\Omega_h\,.
\end{equation} 
Now,  equality \eqref{U.eq3} follows by \eqref{U.eq5}, \eqref{U.eq6}, and \eqref{U.eq7}.  The proof of \eqref{U.eq4} is similar. By a straightforward computation based on the rule of change of variables in integrals  we verify that
\begin{equation}\label{U.eq8}
U^m_h[0,\epsilon_2](\xi)= w_{\Omega^o}[\theta^o[0,\epsilon_2]](0)-w^-_{\tilde\Omega(\epsilon_2)}[\tilde\theta](p^h+\epsilon_2\xi)\qquad\forall\xi\in\mathrm{cl}\Omega_m\,,
\end{equation}
where $\tilde\theta\in C^{1,\alpha}(\partial\tilde\Omega(\epsilon_2))$ is defined by
\[
\tilde\theta(x)\equiv
\theta_h\Big(\frac{x-p^h}{\epsilon_2}\Big)\qquad\forall h\in\{1,2\}\,,\;  x\in\partial\Omega_h(1,\epsilon_2)\,.
\]
By Proposition \ref{L0e=0}, we have 
\begin{equation}\label{U.eq9}
w_{\Omega^o}[\theta^o[0,\epsilon_2]](0)=u^o(0)\,.
\end{equation} 
By Proposition \ref{L0e=0}, by the jump formulae \eqref{jump}, by equality \eqref{HxOe.eq3}, and by definition \eqref{HxOe.eq1}, we have
\[
\begin{split}
-w^-_{\tilde\Omega(\epsilon_2)}[\tilde\theta]_{|\partial\tilde\Omega(\epsilon_2)}&=\tilde f-\sum_{h=1} ^2\biggl(\int_{\partial\tilde\Omega(\epsilon_2)}\tilde f\,\tilde\rho_h\,d\sigma\biggr)\mathcal{X}_{\tilde\Omega(\epsilon_2),h}\\
&=\tilde f-\lim_{y\to\infty}\tilde u(y)\\
&\quad +(H^{1,1}_{\tilde\Omega(\epsilon_2)}-H^{1,2}_{\tilde\Omega(\epsilon_2)})\int_{\partial\Omega_2(1,\epsilon_2)}\nu_{\Omega_2(1,\epsilon_2)}\cdot\nabla\tilde u\,d\sigma\,\mathcal{X}_{\tilde\Omega(\epsilon_2),1}\\
&\quad+(H^{2,2}_{\tilde\Omega(\epsilon_2)}-H^{2,1}_{\tilde\Omega(\epsilon_2)})\int_{\partial\Omega_1(1,\epsilon_2)}\nu_{\Omega_1(1,\epsilon_2)}\cdot\nabla\tilde u\,d\sigma\,\mathcal{X}_{\tilde\Omega(\epsilon_2),2}\,.
\end{split}
\]
Then, by the uniqueness of the solution of the exterior Dirichlet problem, we deduce that 
\begin{equation}\label{U.eq10}
-w^-_{\tilde\Omega(\epsilon_2)}[\tilde\theta]=\tilde u+\tilde w-\lim_{y\to\infty}\tilde u(y)\,.
\end{equation}
Hence, the validity of  \eqref{U.eq4}  follows by \eqref{U.eq8}, \eqref{U.eq9}, and \eqref{U.eq10}.
\qed

\vspace{\baselineskip}

\begin{prop}\label{V}
For all $(\epsilon_1,\epsilon_2)\in ]-\delta_1,\delta_1[\times]-\delta_2,\delta_2[$ we denote by $V[\epsilon_1,\epsilon_2]\equiv(V_1[\epsilon_1,\epsilon_2],V_2[\epsilon_1,\epsilon_2])$ the function of $C^{1,\alpha}(\mathrm{cl}\Omega(\epsilon_1,\epsilon_2))^2$ defined by
\[
V_j[\epsilon_1,\epsilon_2](x)\equiv v_{\Omega^o}[\rho^o_j[\epsilon_1,\epsilon_2]](x)+\sum_{k=1}^2\int_{\partial\Omega_k}S(x-\epsilon_1p^k-\epsilon_1\epsilon_2\eta)\rho_{j,k}[\epsilon_1,\epsilon_2](\eta)\,d\sigma_\eta\quad\forall x\in\mathrm{cl}\Omega(\epsilon_1,\epsilon_2)
\] 
for all $j\in\{1,2\}$.
Then the following statements hold.
\begin{enumerate}
\item[(i)] Let $\Omega_M$ be an open subset of $\Omega^o$ such that $0\notin\mathrm{cl}\Omega_M$. Let $\delta_M\in]0,\delta_1]$ be such that $\mathrm{cl}\Omega_M\cap\mathrm{cl}\Omega_k(\epsilon_1,\epsilon_2)=\emptyset$ for all $(\epsilon_1,\epsilon_2)\in]-\delta_M,\delta_M[\times]-\delta_2,\delta_2[$ and for all $k\in\{1,2\}$. Then there exists a real analytic map $V^M\equiv (V^M_1,V^M_2)$ from $]-\delta_M,\delta_M[\times]-\delta_2,\delta_2[$ to $C^{1,\alpha}(\mathrm{cl}\Omega_M)^2$ such that 
\[
V[\epsilon_1,\epsilon_2](x)=V^M[\epsilon_1,\epsilon_2](x)\qquad\forall x\in\mathrm{cl}\Omega_M\,,\;(\epsilon_1,\epsilon_2)\in]-\delta_M,\delta_M[\times]-\delta_2,\delta_2[\,.
\]
Moreover,
\begin{equation}\label{V.eq1}
V_j^M[\epsilon_1,0](x)=S(x-\epsilon_1p^j)-H^{\Omega^o}_x(\epsilon_1p^j)\quad\forall j\in\{1,2\}\,,\; x\in\mathrm{cl}\Omega_M\,,\;\epsilon_1\in]-\delta_M,\delta_M[\,,
\end{equation}
and 
\begin{equation}\label{V.eq2}
V_j^M[0,\epsilon_2](x)=S(x)-H^{\Omega^o}_x(0)\quad\forall j\in\{1,2\}\,,\; x\in\mathrm{cl}\Omega_M\,,\;\epsilon_2\in]-\delta_2,\delta_2[\,.
\end{equation}
\item[(ii)] Let $h,k\in\{1,2\}$ and $h\neq k$. Let $\Omega_m$ be an open bounded subset of $\mathbb{R}^2\setminus\mathrm{cl}\Omega_h$. Let $\delta_m\in]0,\delta_1]$ be such that $\epsilon_1p^h+\epsilon_1\epsilon_2\mathrm{cl}\Omega_m\subseteq\Omega^o$ and $(\epsilon_1p^h+\epsilon_1\epsilon_2\mathrm{cl}\Omega_m)\cap\mathrm{cl}\Omega_k(\epsilon_1,\epsilon_2)=\emptyset$ for all $(\epsilon_1,\epsilon_2)\in]-\delta_m,\delta_m[^2$. Then there exists a real analytic map $V_h^m\equiv(V^m_{h,1},V^m_{h,2})$ from  $]-\delta_m,\delta_m[^2$ to $C^{1,\alpha}(\mathrm{cl}\Omega_m)$ such that 
\begin{equation}\label{V.eq3}
V_j[\epsilon_1,\epsilon_2](\epsilon_1p^h+\epsilon_1\epsilon_2\xi)=V_{h,j}^m[\epsilon_1,\epsilon_2](\xi)+\delta_{j,h}\frac{1}{2\pi}\log|\epsilon_1\epsilon_2|+\delta_{j,k}\frac{1}{2\pi}\log|\epsilon_1|\quad\forall \xi\in\mathrm{cl}\Omega_m
\end{equation}
for all $j\in\{1,2\}$, $(\epsilon_1,\epsilon_2)\in\left(]-\delta_m,\delta_m[\setminus\{0\}\right)^2$. Moreover,
\begin{equation}\label{V.eq4}
V_{h,j}^m[\epsilon_1,0](\xi)=-H^{\Omega^o}_{\epsilon_1p^h}(\epsilon_1p^j)+\delta_{j,h}\lim_{\eta\to\infty}H^\xi_{\Omega_h}(\eta)+\delta_{j,k}S(p^h-p^k)
\end{equation}
for all $j\in\{1,2\}$, $\xi\in\mathrm{cl}\Omega_m$, and $\epsilon_1\in]-\delta_m,\delta_m[$, and 
\begin{equation}\label{V.eq5}
\begin{split}
&V_{h,j}^m[0,\epsilon_2](\xi) =-H^{\Omega^o}_{0}(0)+\lim_{y\to\infty}H_{\tilde\Omega(\epsilon_2)}^{p^h+\epsilon_2\xi}(y)\\
&\qquad +\left(H^{j,k}_{\tilde\Omega(\epsilon_2)}-H^{j,h}_{\tilde\Omega(\epsilon_2)}\right)\int_{\partial\Omega_h(1,\epsilon_2)}\nu_{\Omega_h(1,\epsilon_2)}(y)\cdot\nabla_y H_{\tilde\Omega(\epsilon_2)}^{p^h+\epsilon_2\xi}(y)\,d\sigma_y-\frac{\log|\epsilon_2|}{2\pi}\delta_{j,h}
\end{split}
\end{equation}
for all $j\in\{1,2\}$, $\xi\in\mathrm{cl}\Omega_m$, and $\epsilon_2\in]-\delta_m,\delta_m[\setminus\{0\}$.
\end{enumerate}
\end{prop}
\proof
To prove statement (i) we take 
\[
V^M[\epsilon_1,\epsilon_2]\equiv V[\epsilon_1,\epsilon_2]_{|\mathrm{cl}\Omega_M}\qquad\forall (\epsilon_1,\epsilon_2)\in]-\delta_M,\delta_M[\times]-\delta_2,\delta_2[\,.
\]
Then, the real analyticity of $V^M$ follows by the standard properties of integral  operators with real analytic kernels and with no singularity (see, e.g., Lanza de Cristoforis and the second author \cite{LaMu13}), by the classical mapping properties of layer potentials (cf., {\it e.g.}, Miranda \cite{Mi65}), and by Proposition \ref{rhoee}. The validity of equalities \eqref{V.eq1} and \eqref{V.eq2} can be deduced by Proposition \ref{HOox} and by Propositions \ref{Me0=0} and \ref{M0e=0}.

We now consider statement (ii). We define 
\[
\begin{split}
V^m_{h,j}[\epsilon_1,\epsilon_2](\xi)&\equiv v_{\Omega^o}[\rho^o_j[\epsilon_1,\epsilon_2]](\epsilon_1p^h+\epsilon_1\epsilon_2\xi)+v_{\Omega_h}[\rho_{j,h}[\epsilon_1,\epsilon_2]](\xi)\\
&\quad +\int_{\partial\Omega_k}S(p^h-p^k-\epsilon_2(\xi-\eta))\rho_{j,k}[\epsilon_1,\epsilon_2](\eta)\,d\sigma_\eta\qquad\forall\xi\in\mathrm{cl}\Omega_m
\end{split}
\]
for all $j\in\{1,2\}$ and $(\epsilon_1,\epsilon_2)\in]-\delta_m,\delta_m[^2$. Then, by the standard properties of integral  operators with real analytic kernels and with no singularity (see, e.g., Lanza de Cristoforis and the second author \cite{LaMu13}), by the  mapping properties of layer potentials (cf., {\it e.g.}, Miranda \cite{Mi65}), and by Proposition \ref{rhoee} we verify that $V^m_h\equiv(V^m_{h,1},V^m_{h,2})$ is real analytic. Then  equality \eqref{V.eq3} follows by a straightforward computation based on the rule of change of variables in integrals and on Proposition \ref{Mee=0}. 
To prove equality \eqref{V.eq4} we observe that by Proposition \ref{Me0=0}
\[
V^m_{h,j}[\epsilon_1,0](\xi)= v_{\Omega^o}[\rho^o_j[\epsilon_1,0]](\epsilon_1p^h)+v_{\Omega_h}[\rho_{j,h}[\epsilon_1,0]](\xi)+S(p^h-p^k)\delta_{j,k}\quad\forall\xi\in\mathrm{cl}\Omega_m\,.
\]
Then the validity of \eqref{V.eq4} follows by Proposition \ref{HOox} and equality \eqref{HxOh.eq0}. By Proposition \ref{HOox} and by equality \eqref{HxOe.eq2} one verifies \eqref{V.eq5}.
\qed

\vspace{\baselineskip}

\begin{prop}\label{F}
Let $F\equiv(F_1,F_2)$ be the function from $]-\delta_1,\delta_1[\times]-\delta_2,\delta_2[$ to $\mathbb{R}^2$ defined by
\[
F_j[\epsilon_1,\epsilon_2]\equiv \int_{\partial\Omega^o} f^o\rho^o_j[\epsilon_1,\epsilon_2]\,d\sigma+\sum_{h=1}^2\int_{\partial\Omega_h}f_h\, \rho_{j,h}[\epsilon_1,\epsilon_2]\,d\sigma\qquad\forall j\in\{1,2\}\,,\;(\epsilon_1,\epsilon_2)\in]-\delta_1,\delta_1[\times]-\delta_2,\delta_2[\,.
\]
Then $F$ is real analytic. Moreover,  we have
\begin{align}\label{F.eq1}
&F_j[\epsilon_1,0]=-u^o(\epsilon_1p^j)+\lim_{y\to\infty} u_j(y)\,,\\
\label{F.eq2}
&F_j[0,\epsilon_2]=-u^o(0)+\lim_{y\to\infty} \tilde u(y)+(H^{j,k}_{\tilde\Omega(\epsilon_2)}-H^{j,h}_{\tilde\Omega(\epsilon_2)})\int_{\partial\Omega_k(1,\epsilon_2)}\nu_{\Omega_k(1,\epsilon_2)}\cdot\nabla\tilde u\,d\sigma\,,
\end{align}
for all $j,h,k\in\{1,2\}$, $h\neq k$, $\epsilon_1\in]-\delta_1,\delta_1[$, and $\epsilon_2\in]-\delta_2,\delta_2[\setminus\{0\}$.
\end{prop}
\proof 
The real analyticity of $F$ is a consequence of Proposition \ref{rhoee}. The validity of \eqref{F.eq1} follows by \eqref{HOox.eq1} and \eqref{HxOh.eq1}. To prove \eqref{F.eq2} one observes that 
\[
F_j[0,\epsilon_2]=\int_{\partial\Omega^o} f^o\rho^o_j[0,\epsilon_2]\,d\sigma+\int_{\partial\tilde\Omega(\epsilon_2)}\tilde f\, \tilde\rho_{j}\,d\sigma
\] with $\tilde{f}(x)\equiv
f_h((x-p^h)/\epsilon_2)$ for all $h\in\{1,2\}$ and $x\in\partial\Omega_h(1,\epsilon_2)$ and $\tilde\rho_{j}$ as in Proposition \ref{M0e=0}. Then the validity of  \eqref{F.eq2} follows by \eqref{HOox.eq2},   \eqref{HxOe.eq1}, and  \eqref{HxOe.eq3}.
\qed

\vspace{\baselineskip}

Here below  $M_{2\times 2}(\mathbb{R})$ denotes the space of the $2\times 2$ real matrices.

\begin{prop}\label{R}
Let $R\equiv(R_{i,j})_{(i,j)\in\{1,2\}^2}$ be the function from $]-\delta_1,\delta_1[\times]-\delta_2,\delta_2[$ to $M_{2\times 2}(\mathbb{R})$   defined by
\[
\begin{split}
R_{i,j}[\epsilon_1,\epsilon_2]&\equiv\frac{1}{|\partial\Omega_j|}\int_{\partial\Omega_j}v_{\Omega^o}[\rho^o_i[\epsilon_1,\epsilon_2]](\epsilon_1p^j+\epsilon_1\epsilon_2\xi)+v_{\Omega_h}[\rho_{i,j}[\epsilon_1,\epsilon_2]](\xi)\\
&\qquad\qquad\quad+\left(\int_{\partial\Omega_k}S(p^j-p^k+\epsilon_2(\xi-\eta))\rho_{i,k}[\epsilon_1,\epsilon_2](\eta)\,d\sigma_\eta\right)\, d\sigma_\xi
\end{split}
\]
for all $(\epsilon_1,\epsilon_2)\in]-\delta_1,\delta_1[\times]-\delta_2,\delta_2[$ and for all $i,j,k\in\{1,2\}$ with $j\neq k$. Then $R$ is real analytic and 
\begin{align}
\label{R.eq1}
R_{i,j}[\epsilon_1,0]&=-H^{\Omega^o}_{\epsilon_1p^j}(\epsilon_1p^i)+\delta_{i,j}\lim_{\eta\to\infty} H_{\Omega_i}^0(\eta)+(1-\delta_{i,j})S(p^i-p^j)\,,\\
\label{R.eq2}
R_{i,j}[0,\epsilon_2]&=-H^{\Omega^o}_0(0)+H^{i,j}_{\tilde\Omega(\epsilon_2)}-\frac{\log|\epsilon_2|}{2\pi}\delta_{i,j}
\end{align}
for all $i,j\in\{1,2\}$, $\epsilon_1\in]-\delta_1,\delta_1[$, and $\epsilon_2\in]-\delta_2,\delta_2[\setminus\{0\}$.
\end{prop}
\proof
The real analyticity of $R$ is a consequence of Proposition \ref{rhoee} and of the mapping properties of the single layer potential. Equality \eqref{R.eq1} follows by Proposition \ref{HOox}, by \eqref{HxOh.eq0}, and by Proposition \ref{Mee=0}. Equality \eqref{R.eq2} follows by  Proposition \ref{HOox} and by equality \eqref{HxOe.eq1.1}.
\qed

\vspace{\baselineskip}

\begin{prop}\label{Lambda}
Let $\epsilon_1\in]-\delta_1,\delta_1[\setminus\{0\}$ and $\epsilon_2\in]-\delta_2,\delta_2[\setminus\{0\}$. Then  the matrix $\Lambda(\epsilon_1,\epsilon_2)\equiv \bigl(\Lambda_{i,j}(\epsilon_1,\epsilon_2)\bigr)_{(i,j)\in\{1,2\}^2}$ defined by
\[
\Lambda_{i,j}(\epsilon_1,\epsilon_2)\equiv \delta_{i,j}\frac{1}{2\pi}\log|\epsilon_1\epsilon_2|+(1-\delta_{i,j})\frac{1}{2\pi}\log|\epsilon_1|+R_{i,j}[\epsilon_1,\epsilon_2]\quad\forall i,j\in\{1,2\}
\]
satisfies the equality 
\begin{equation}\label{Lambda.eq1}
\Lambda_{i,j}(\epsilon_1,\epsilon_2)=\frac{1}{|\partial\Omega_j(\epsilon_1,\epsilon_2)|}\int_{\partial\Omega_j(\epsilon_1,\epsilon_2)}v_{\Omega(\epsilon_1,\epsilon_2)}[\tau_i]\,d\sigma\qquad\qquad\forall i,j\in\{1,2\}
\end{equation}
with $\tau_i\in C^{0,\alpha}(\partial\Omega(\epsilon_1,\epsilon_2))$ defined as in \eqref{taui}. In particular, the matrix $\Lambda(\epsilon_1,\epsilon_2)$ is invertible.
\end{prop}
\proof Equality \eqref{Lambda.eq1} follows by Proposition \ref{Mee=0} and by the rule of change of variables in integrals.  The invertibility of $\Lambda(\epsilon_1,\epsilon_2)$ is a consequence of Lemma \ref{LambdaO}. \qed

\vspace{\baselineskip}

We are now ready to prove our main  Theorem \ref{ue1e2}, where we introduce representation formulas for $u_{\epsilon_1,\epsilon_2}$ and for $u_{\epsilon_1,\epsilon_2}(\epsilon_1p^h+\epsilon_1\epsilon_2\,\cdot\,)$ in terms of real analytic functions of the pair $(\epsilon_1,\epsilon_2)$ and of  elementary functions of $\log|\epsilon_1|$ and $\log|\epsilon_1\epsilon_2|$. In the sequel,  $A^t$ denotes the transpose of a matrix $A$ and $A^{-1}$ denotes the inverse of an invertible matrix $A$.

\begin{thm}\label{ue1e2}
The following statements hold.
\begin{enumerate}
\item[(i)]  Let $\Omega_M$ be an open subset of $\Omega^o$ such that $0\notin\mathrm{cl}\Omega_M$. Let $\delta_M\in]0,\delta_1]$ be such that $\mathrm{cl}\Omega_M\cap\mathrm{cl}\Omega_k(\epsilon_1,\epsilon_2)=\emptyset$ for all $(\epsilon_1,\epsilon_2)\in]-\delta_M,\delta_M[\times]-\delta_2,\delta_2[$ and for all $k\in\{1,2\}$. Then 
\[
u_{\epsilon_1,\epsilon_2|\mathrm{cl}\Omega_M}=u^o+\epsilon_1\epsilon_2\;U^M[\epsilon_1,\epsilon_2]+F[\epsilon_1,\epsilon_2]^t\,\Lambda(\epsilon_1,\epsilon_2)^{-1}\, V^M[\epsilon_1,\epsilon_2]
\]
 for all $\epsilon_1\in]-\delta_M,\delta_M[\setminus\{0\}$ and $\epsilon_2\in]-\delta_2,\delta_2[\setminus\{0\}$.
\item[(ii)]  Let $h,k\in\{1,2\}$ and $h\neq k$. Let $\Omega_m$ be an open bounded subset of $\mathbb{R}^2\setminus\mathrm{cl}\Omega_h$. Let $\delta_m\in]0,\delta_1]$ be such that $\epsilon_1p^h+\epsilon_1\epsilon_2\mathrm{cl}\Omega_m\subseteq\Omega^o$ and $(\epsilon_1p^h+\epsilon_1\epsilon_2\mathrm{cl}\Omega_m)\cap\mathrm{cl}\Omega_k(\epsilon_1,\epsilon_2)=\emptyset$ for all $(\epsilon_1,\epsilon_2)\in]-\delta_m,\delta_m[^2$. Then  
\[
u_{\epsilon_1,\epsilon_2}(\epsilon_1p^h+\epsilon_1\epsilon_2\,\cdot\,)_{|\mathrm{cl}\Omega_m}=U_h^m[\epsilon_1,\epsilon_2]+F[\epsilon_1,\epsilon_2]^t\, \Lambda(\epsilon_1,\epsilon_2)^{-1} \, \left(V_h^m[\epsilon_1,\epsilon_2]+\mathcal{S}_h(\epsilon_1,\epsilon_2)\right)
\]
 for all $(\epsilon_1,\epsilon_2)\in(]-\delta_m,\delta_m[\setminus\{0\})^2$, where $\mathcal{S}_h(\epsilon_1,\epsilon_2)\in\mathbb{R}^2$ is defined by
 \[
\mathcal{S}_h(\epsilon_1,\epsilon_2)_j\equiv \delta_{j,h}\frac{1}{2\pi}\log|\epsilon_1\epsilon_2|+(1-\delta_{j,h})\frac{1}{2\pi}
\log|\epsilon_1|\qquad\forall j\in\{1,2\}\,.
 \]
 \end{enumerate}
\end{thm}
\proof 
By Propositions \ref{ugO} we have
 \[
u_{\epsilon_1,\epsilon_2}(x)\equiv w_{\Omega(\epsilon_1,\epsilon_2)}^+[\phi](x)+\sum_{i,j=1}^2\biggl(\int_{\partial\Omega(\epsilon_1,\epsilon_2)}f\tau_i\,d\sigma\biggr)(\Lambda(\epsilon_1,\epsilon_2)^{-1})_{i,j}\, v_{\Omega(\epsilon_1,\epsilon_2)}[\tau_j](x)\quad\forall x\in\mathrm{cl}{\Omega(\epsilon_1,\epsilon_2)}
\]
with $\phi$ as in Proposition \ref{Lee=0},  $\tau_1$ and $\tau_2$ as in \eqref{taui},  and $\Lambda(\epsilon_1,\epsilon_2)$ as in \eqref{Lambda.eq1}. Then the validity of (i) and (ii) follows by Propositions \ref{U}, \ref{V},  \ref{F}, and \ref{Lambda}, and by a computation based on the rule of change of variables in integrals.
\qed

\section{Asymptotic behaviour of $u_{\epsilon_1,\epsilon_2}$ as $(\epsilon_1,\epsilon_2)\to(0,\gamma_0)$}\label{as}

In this section we show how Theorem \ref{ue1e2} can be exploited to obtain asymptotic approximations of the solution of problem \eqref{dir} as the pair of parameters $(\epsilon_1,\epsilon_2)$ approaches  a degenerate pair $(0,\gamma_0)$. As we shall see, the function $1/\log|\epsilon_1\epsilon_2|$ will appear in many of our expressions and, in order that such expressions make sense,  we have to ensure that $|\epsilon_1\epsilon_2|<1$ in the admissible set. Then, we shrink $\delta_1$ and we  assume that in this section we have 
\[
\delta_1\in]0,1/\delta_2[\,.
\]
In the following Proposition \ref{Lambda-1} we describe the inverse matrix $\Lambda(\epsilon_1,\epsilon_2)^{-1}$. In the sequel,  $A^*$ denotes the cofactor matrix of a matrix $A$, so that $A^{*t}$ is the adjugate of  $A$. 

\begin{prop}\label{Lambda-1}
Let $\epsilon_1\in]-\delta_1,\delta_1[\setminus\{0\}$ and $\epsilon_2\in]-\delta_2,\delta_2[\setminus\{0\}$. Then we have
\[
\det\Lambda(\epsilon_1,\epsilon_2)=\frac{1}{4\pi^2}\;\mathcal{R}_{\epsilon_1,\epsilon_2}\;\log|\epsilon_1\epsilon_2|
\]
and
\[
\Lambda(\epsilon_1,\epsilon_2)^{-1}=\frac{2\pi}{\mathcal{R}_{\epsilon_1,\epsilon_2}}
\left\{\left(
\begin{array}{cc}
1&-\frac{\log|\epsilon_1|}{\log|\epsilon_1\epsilon_2|}\\
-\frac{\log|\epsilon_1|}{\log|\epsilon_1\epsilon_2|}&1
\end{array}
\right)
+2\pi\frac{1}{\log|\epsilon_1\epsilon_2|}R[\epsilon_1,\epsilon_2]^{*t}\right\}
\]
with
\begin{equation}\label{Lambda-1.eq1}
\begin{split}
\mathcal{R}_{\epsilon_1,\epsilon_2}\equiv& \log|\epsilon_2|+\frac{\log|\epsilon_1|}{\log|\epsilon_1\epsilon_2|}\log|\epsilon_2|\\
&-2\pi (R_{1,2}[\epsilon_1,\epsilon_2]+R_{2,1}[\epsilon_1,\epsilon_2])\frac{\log|\epsilon_1|}{\log|\epsilon_1\epsilon_2|}\\
&+2\pi (R_{1,1}[\epsilon_1,\epsilon_2]+R_{2,2}[\epsilon_1,\epsilon_2])\\
&+4\pi^2 (R_{1,1}[\epsilon_1,\epsilon_2]R_{2,2}[\epsilon_1,\epsilon_2]-R_{1,2}[\epsilon_1,\epsilon_2]R_{2,1}[\epsilon_1,\epsilon_2])\frac{1}{\log|\epsilon_1\epsilon_2|}\,.
\end{split}
\end{equation}
\end{prop}

We observe that, since $\Lambda(\epsilon_1,\epsilon_2)$ is invertible by Proposition \ref{Lambda}, we have that $\mathcal{R}_{\epsilon_1,\epsilon_2}\neq 0$ for all  $\epsilon_1\in]-\delta_1,\delta_1[\setminus\{0\}$ and $\epsilon_2\in]-\delta_2,\delta_2[\setminus\{0\}$.

In the following Proposition \ref{tildeV} we write a convenient expression for $u_{\epsilon_1,\epsilon_2|\mathrm{cl}\Omega_M}$. We exploit the following definition
\begin{equation}\label{Green}
G^{\Omega^o}(x,y)\equiv S(x-y)-H^{\Omega^o}_x(y)\qquad\forall x,y\in\Omega^o\,,\; x\neq y
\end{equation}
(cf.~Proposition \ref{HOox}). We observe that $G^{\Omega^o}$ is the Dirichlet Green function for the domain $\Omega^o$.

\begin{prop} \label{tildeV}
Let $\Omega_M$ be an open subset of $\Omega^o$ such that $0\notin\mathrm{cl}\Omega_M$. Let $\delta_M\in]0,\delta_1]$ be such that $\mathrm{cl}\Omega_M\cap\mathrm{cl}\Omega_k(\epsilon_1,\epsilon_2)=\emptyset$ for all $(\epsilon_1,\epsilon_2)\in]-\delta_M,\delta_M[\times]-\delta_2,\delta_2[$ and for all $k\in\{1,2\}$. Then there exists a real analytic map ${X}^M\equiv (X^M_1,X^M_2)$ from $]-\delta_M,\delta_M[\times]-\delta_2,\delta_2[$ to $C^{1,\alpha}(\mathrm{cl}\Omega_M)^2$ such that
\begin{equation}\label{tildeV.eq1}
\begin{split}
u_{\epsilon_1,\epsilon_2}(x)&=u^o(x)+\epsilon_1\epsilon_2\;U^M[\epsilon_1,\epsilon_2](x)\\
&\quad+2\pi\frac{1}{\log|\epsilon_1\epsilon_2|}\frac{\log|\epsilon_2|}{\mathcal{R}_{\epsilon_1,\epsilon_2}}\left(F_1[\epsilon_1,\epsilon_2]+F_2[\epsilon_1,\epsilon_2]\right)
G^{\Omega^o}(x,0)\\
&\quad
+2\pi
\frac{\epsilon_1}{\mathcal{R}_{\epsilon_1,\epsilon_2}}
F[\epsilon_1,\epsilon_2]^t
\left(
\begin{array}{cc}
1&-\frac{\log|\epsilon_1|}{\log|\epsilon_1\epsilon_2|}\\
-\frac{\log|\epsilon_1|}{\log|\epsilon_1\epsilon_2|}&1
\end{array}
\right)X^M[\epsilon_1,\epsilon_2](x)
\\
&\quad
+4\pi^2
\frac{1}{\log|\epsilon_1\epsilon_2|}
\frac{1}{\mathcal{R}_{\epsilon_1,\epsilon_2}}
F[\epsilon_1,\epsilon_2]^t\, R[\epsilon_1,\epsilon_2]^{*t}\, V^M[\epsilon_1,\epsilon_2](x)\qquad\forall x\in\mathrm{cl}\Omega^M
\end{split}
\end{equation}
for all $\epsilon_1\in ]-\delta_M,\delta_M[\setminus\{0\}$ and $\epsilon_2\in]-\delta_2,\delta_2[\setminus\{0\}$.
\end{prop} 
\proof 
By Proposition \ref{V} (i) and by standard properties of real analytic functions there exists a real analytic map $X^M\equiv (X^M_1,X^M_2)$ from $]-\delta_M,\delta_M[\times]-\delta_2,\delta_2[$ to $C^{1,\alpha}(\mathrm{cl}\Omega_M)^2$ such that  
\[
V^M_j[\epsilon_1,\epsilon_2](x)= G^{\Omega^o}(x,0)+\epsilon_1X^M_j[\epsilon_1,\epsilon_2](x)\quad\forall x\in\mathrm{cl}\Omega^M\,,\; (\epsilon_1,\epsilon_2)\in ]-\delta_M,\delta_M[\times]-\delta_2,\delta_2[\,,\; j\in\{1,2\}\,.
\]
Then, by a straightforward computation based on Proposition \ref{Lambda-1} we have
\[
\begin{split}
&\Lambda(\epsilon_1,\epsilon_2)^{-1}V^M[\epsilon_1,\epsilon_2](x)\\
&\qquad=2\pi\frac{1}{\log|\epsilon_1\epsilon_2|}\frac{\log|\epsilon_2|}{\mathcal{R}_{\epsilon_1,\epsilon_2}}\left(
\begin{array}{c}
G^{\Omega^o}(x,0)\\
G^{\Omega^o}(x,0)
\end{array}
\right)\\
&\qquad\quad
+2\pi\frac{\epsilon_1}{\mathcal{R}_{\epsilon_1,\epsilon_2}}
\left(
\begin{array}{cc}
1&-\frac{\log|\epsilon_1|}{\log|\epsilon_1\epsilon_2|}\\
-\frac{\log|\epsilon_1|}{\log|\epsilon_1\epsilon_2|}&1
\end{array}
\right)X^M[\epsilon_1,\epsilon_2](x)
\\
&\qquad\quad
+4\pi^2
\frac{1}{\log|\epsilon_1\epsilon_2|}
\frac{1}{\mathcal{R}_{\epsilon_1,\epsilon_2}}
R[\epsilon_1,\epsilon_2]^{*t}\, V^M[\epsilon_1,\epsilon_2](x)\qquad\forall x\in\mathrm{cl}\Omega^M
\end{split}
\]
for all $\epsilon_1\in ]-\delta_M,\delta_M[\setminus\{0\}$ and $\epsilon_2\in]-\delta_2,\delta_2[\setminus\{0\}$. Now the validity of the statement follows by Theorem \ref{ue1e2}.
\qed

\vspace{\baselineskip}

We now observe that if  we try to pass to the limit in the representation formula \eqref{tildeV.eq1} we face the problem that
 \[
 \lim_{(\epsilon_1,\epsilon_2) \to (0,\gamma_0)}\frac{\log |\epsilon_1|}{\log|\epsilon_1\epsilon_2|}
 \]
does not exist when $\gamma_0=0$. As it has been announced in the introduction,  we can overcome this difficulty by replacing $\epsilon_1$ with a positive parameter $t$ and by taking $\epsilon_2=\gamma(t)$, where $\gamma$ is a function from a right neighbourhood of $0$ to $]0,\delta_2[$ such that the limits
\begin{equation}\label{limits}
\gamma_0\equiv\lim_{t\to 0}\gamma(t)\qquad\text{and}\qquad\lambda_0\equiv\lim_{t\to 0}\frac{\log t}{\log(t\gamma(t))}
\end{equation}
exist finite in $[0,\delta_2[$ and $[0,+\infty[$, respectively. Then we investigate the first and second term in the asymptotic expansion of $u_{t,\gamma(t)}$ as $t\to 0^+$. We observe that we have to distinguish the case when $\lim_{t\to 0^+}\gamma(t)=0$ from the case when $\lim_{t\to 0^+}\gamma(t)>0$.  We shall also need the following technical lemma. 

\begin{lem}\label{cgamma0}
Let $\gamma_0\in]0,\delta_2[$. Let $c_{\gamma_0}\in \mathbb{R}$ be defined by
\[
c_{\gamma_0}\equiv H^{1,1}_{\tilde\Omega(\gamma_0)}-H^{1,2}_{\tilde\Omega(\gamma_0)}-H^{2,1}_{\tilde\Omega(\gamma_0)}+H^{2,2}_{\tilde\Omega(\gamma_0)}\,.
\]
Then $c_{\gamma_0}\neq 0$.
\end{lem}
\proof
By \eqref{HxOe.eq1} we have 
\[
c_{\gamma_0}=v_{\tilde\Omega(\gamma_0)}[\tilde\rho_1-\tilde\rho_2](p^1)-v_{\tilde\Omega(\gamma_0)}[\tilde\rho_1-\tilde\rho_2](p^2)\,.
\]
where $\tilde\rho_1, \tilde\rho_2\in C^{0,\alpha}(\partial\tilde\Omega(\epsilon_2))$ are defined as in \eqref{tilderhoi}. By Proposition \ref{M0e=0}, $\tilde\rho_1-\tilde\rho_2$ belongs to $\mathrm{Ker}(-\frac{1}{2}I_{\tilde\Omega(\gamma_0)}+W^*_{\tilde\Omega(\gamma_0)})$. Then, the jump formula for $v^+_{\tilde\Omega(\gamma_0)}[\tilde\rho_1-\tilde\rho_2]$ in \eqref{jump} implies that  $v_{\tilde\Omega(\gamma_0)}[\tilde\rho_1-\tilde\rho_2]$ is constant on $\mathrm{cl}\Omega_1(0,\gamma_0)$ and  on $\mathrm{cl}\Omega_2(0,\gamma_0)$. It follows that $c_{\gamma_0}=0$ only if $v_{\tilde\Omega(\gamma_0)}[\tilde\rho_1-\tilde\rho_2]$ equals the same constant on $\mathrm{cl}\Omega_1(0,\gamma_0)$ and  on $\mathrm{cl}\Omega_2(0,\gamma_0)$. That is, if  $v_{\tilde\Omega(\gamma_0)}[\tilde\rho_1-\tilde\rho_2]$ is constant on the whole of $\mathrm{cl}\tilde\Omega(\gamma_0)$.  Then we observe that by Proposition \ref{M0e=0} we also have
\begin{equation}\label{cgamma0.eq1}
\int_{\partial\Omega_i(1,\gamma_0)}\tilde\rho_1-\tilde\rho_2\, d\sigma=(-1)^{i+1}\,.
\end{equation}
Thus $\int_{\tilde\Omega(\gamma_0)}\tilde\rho_1-\tilde\rho_2\, d\sigma=0$ and $\tilde\rho_1-\tilde\rho_2\in \bigl(\mathrm{Ker}(-\frac{1}{2}I_{\tilde\Omega(\gamma_0)}+W^*_{\tilde\Omega(\gamma_0)})\bigr)_0$. So, by Lemma \ref{ker} (vii),  we deduce that $c_{\gamma_0}= 0$ only for $\tilde\rho_1=\tilde\rho_2$. However, the latter equality is in contradiction with \eqref{cgamma0.eq1}. Thus $c_{\gamma_0}\neq 0$. 
\qed

\vspace{\baselineskip}

We now prove our main result on the asymptotic behaviour of $u_{t,\gamma(t)}$ as $t\to 0^+$.

\begin{prop}\label{gamma}
Let $\Omega_M$ be an open subset of $\Omega^o$ such that $0\notin\mathrm{cl}\Omega_M$. Let $\delta_M\in]0,\delta_1]$ be such that $\mathrm{cl}\Omega_M\cap\mathrm{cl}\Omega_k(\epsilon_1,\epsilon_2)=\emptyset$ for all $(\epsilon_1,\epsilon_2)\in]-\delta_M,\delta_M[\times]-\delta_2,\delta_2[$ and for all $k\in\{1,2\}$. Let $\delta^*_M\in]0,\delta_M]$. Let $\gamma$ be a function from $]0,\delta^*_M[$ to $]0,\delta_2[$ such that the limits in \eqref{limits}
exist finite in $[0,\delta_2[$ and $[0,+\infty[$, respectively. Then the following statements hold:
\begin{enumerate}
\item[(i)] If $\gamma_0=0$, then we have
\[
\begin{split}
&u_{t,\gamma(t)|\mathrm{cl}\Omega_M}=u^o_{|\mathrm{cl}\Omega_M}\\
&\ +\frac{1}{\log(t\gamma(t))}\frac{2\pi}{1+\lambda_0}\Bigl(\lim_{y\to\infty} u_1(y)+\lim_{y\to\infty} u_2(y)-2u^o(0)\Bigr)G^{\Omega^o}(\cdot,0)_{|\mathrm{cl}\Omega_M}+o\left(\frac{1}{\log(t\gamma(t))}\right)
\end{split}
\]
as $t\to 0^+$.
\item[(ii)] If $\gamma_0\in]0,\delta_2[$, then $\lambda_0=1$ and
\begin{equation}\label{gamma.eq2}
\begin{split}
&u_{t,\gamma(t)|\mathrm{cl}\Omega_M}=u^o_{|\mathrm{cl}\Omega_M}\\
 & +\frac{2\pi}{\log t}\biggl(\lim_{y\to\infty} \tilde u(y)-u^o(0)+\biggl.\left(H^{2,1}_{\tilde\Omega(\gamma_0)}-H^{1,2}_{\tilde\Omega(\gamma_0)}\right)\int_{\Omega_1(1,\gamma_0)}\nu_{\Omega_1(1,\gamma_0)}\cdot\nabla\tilde u\, d\sigma\biggr) G^{\Omega^o}(\cdot,0)_{|\mathrm{cl}\Omega_M}+o\left(\frac{1}{\log t}\right)
\end{split}
\end{equation}
 as $t\to 0^+$.
\end{enumerate}
\end{prop}
\proof We first prove (i). If $\gamma_0=0$, then we have 
\[
\lim_{t\to 0^+}\frac{\mathcal{R}_{t,\gamma(t)}}{\log\gamma(t)}=1+\lambda_0
\] 
(cf.~\eqref{Lambda-1.eq1}). Then the validity of (i) follows by Proposition \ref{tildeV}, by the membership of $t\gamma(t)$, $t/\log\gamma(t)$, and $1/(\log(t\gamma(t))\log\gamma(t))$ in $o(1/\log(t\gamma(t)))$, and by a straightforward computation. We now pass to consider (ii). First we observe that the condition $\gamma_0\in]0,\delta_2[$ readily implies that 
$\lambda_0=1$.
Then, by \eqref{Lambda-1.eq1} we deduce that
\[
\lim_{t\to 0^+}\mathcal{R}_{t,\gamma(t)}=2\log\gamma_0+2\pi\bigl(R_{1,1}[0,\gamma_0]+R_{2,2}[0,\gamma_0]-R_{1,2}[0,\gamma_0]-R_{2,1}[0,\gamma_0]\bigr)\,.
\] 
Thus, \eqref{R.eq2} implies that
\begin{equation}\label{gamma.eq3}
\lim_{t\to 0^+}\frac{1}{\mathcal{R}_{t,\gamma(t)}}=\frac{1}{2\pi c_{\gamma_0}}\,.
\end{equation}
Next, by \eqref{F.eq2}  we verify that
\begin{equation}\label{gamma.eq4}
\begin{split}
(F_1[0,\gamma_0]+F_2[0,\gamma_0])
G^{\Omega^o}(x,0) =& \left(2\lim_{y\to\infty} \tilde u(y)-2u^o(0)\right)G^{\Omega^o}(x,0)\\
&+d_{\gamma_0}\left(\int_{\Omega_1(1,\gamma_0)}\nu_{\Omega_1(1,\gamma_0)}\cdot\nabla\tilde u\, d\sigma\right)G^{\Omega^o}(x,0)
\end{split}
\end{equation}
for all $x\in\mathrm{cl}\Omega_M$, with 
$d_{\gamma_0}\equiv\left(H^{1,1}_{\tilde\Omega(\gamma_0)}-H^{1,2}_{\tilde\Omega(\gamma_0)}+H^{2,1}_{\tilde\Omega(\gamma_0)}-H^{2,2}_{\tilde\Omega(\gamma_0)}\right)$.
By \eqref{V.eq2}, \eqref{F.eq2}, \eqref{R.eq2}, and  by equality 
\[
\begin{split}
&c_{\gamma_0}\left(H^{2,1}_{\tilde\Omega(\gamma_0)}-H^{1,2}_{\tilde\Omega(\gamma_0)}\right)\\
&\qquad=\left(H^{1,1}_{\tilde\Omega(\gamma_0)}-H^{1,2}_{\tilde\Omega(\gamma_0)}\right)\left(H^{2,2}_{\tilde\Omega(\gamma_0)}-H^{1,2}_{\tilde\Omega(\gamma_0)}\right)+\left(H^{2,1}_{\tilde\Omega(\gamma_0)}-H^{2,2}_{\tilde\Omega(\gamma_0)}\right)\left(H^{1,1}_{\tilde\Omega(\gamma_0)}-H^{2,1}_{\tilde\Omega(\gamma_0)}\right)
\end{split}
\]
 we compute that
\begin{equation}\label{gamma.eq5}
\begin{split}
F[0,\gamma_0]^t\, &R[0,\gamma_0]^{*t}\,V^M[0,\gamma_0](x)\\
&= c_{\gamma_0} \left(\lim_{y\to\infty} \tilde u(y)-u^o(0)\right)G^{\Omega^o}(x,0)\\
&\quad -\frac{\log\gamma_0}{2\pi}\left(2\lim_{y\to\infty} \tilde u(y)-2u^o(0)\right)G^{\Omega^o}(x,0)\\
&\quad +c_{\gamma_0}\left(H^{2,1}_{\tilde\Omega(\gamma_0)}-H^{1,2}_{\tilde\Omega(\gamma_0)}\right)\left(\int_{\Omega_1(1,\gamma_0)}\nu_{\Omega_1(1,\gamma_0)}\cdot\nabla\tilde u\, d\sigma\right) G^{\Omega^o}(x,0)\\
&\quad -d_{\gamma_0}\frac{\log\gamma_0}{2\pi}\left(\int_{\Omega_1(1,\gamma_0)}\nu_{\Omega_1(1,\gamma_0)}\cdot\nabla\tilde u\, d\sigma\right) G^{\Omega^o}(x,0)
\end{split}
\end{equation}
for all $x\in\mathrm{cl}\Omega_M$.
Now, the validity of \eqref{gamma.eq2} follows by \eqref{tildeV.eq1},  by \eqref{gamma.eq3}--\eqref{gamma.eq5}, and by the asymptotic formula
\[
\frac{1}{\log(t\gamma(t))}=\frac{1}{\log t}+o\left(\frac{1}{\log t}\right)\qquad\text{as }t\to 0^+\,.
\]
\qed

\vspace{\baselineskip}

 We observe that the factor
 \[
 \left(H^{2,1}_{\tilde\Omega(\gamma_0)}-H^{1,2}_{\tilde\Omega(\gamma_0)}\right)\int_{\Omega_1(1,\gamma_0)}\nu_{\Omega_1(1,\gamma_0)}\cdot\nabla\tilde u\, d\sigma
 \] 
 appearing in \eqref{gamma.eq2} vanishes when 
\begin{equation}\label{cond1}
\int_{\Omega_1(1,\gamma_0)}\nu_{\Omega_1(1,\gamma_0)}\cdot\nabla\tilde u\, d\sigma=0\,,
\end{equation}
a condition which is equivalent to $\int_{\Omega_2(1,\gamma_0)}\nu_{\Omega_2(1,\gamma_0)}\cdot\nabla\tilde u\, d\sigma=0$, because
\[
\int_{\Omega_2(1,\gamma_0)}\nu_{\Omega_2(1,\gamma_0)}\cdot\nabla\tilde u\, d\sigma=-\int_{\Omega_1(1,\gamma_0)}\nu_{\Omega_1(1,\gamma_0)}\cdot\nabla\tilde u\, d\sigma\,.
\]
It also vanishes for 
\[
H^{2,1}_{\tilde\Omega(\gamma_0)}=H^{1,2}_{\tilde\Omega(\gamma_0)},
\]
{\it i.e.} for
\begin{equation}\label{cond3}
v_{\tilde\Omega(\gamma_0)}[\tilde\rho_2](p^1)=v_{\tilde\Omega(\gamma_0)}[\tilde\rho_1](p^2)\,.
\end{equation}
Condition \eqref{cond1} concerns $\tilde u$ and thus depends on the geometry of the holes and on the boundary data  $f_1$ and $f_2$.  It is verified for example when $\Omega_2=-\Omega_1$ and $f_2(x)=f_1(-x)$ for all $x\in\partial\Omega_2$.
 Instead, $\tilde\rho_1$ and $\tilde\rho_2$ depend only on the geometry of the holes (see Proposition \ref{M0e=0}). Accordingly, \eqref{cond3} is a geometric conditions on the holes. A simple arguments shows that it is verified for example when $\Omega_2=-\Omega_1$.
 
To conclude, we observe that an analog of Proposition \ref{gamma} can also be proved for the microscopic behaviour of the solution near the boundaries of the holes. Then one can exploit such results to investigate the asymptotic behaviour of functionals of the solution. For example, one may consider the energy integral, which plays an important role in the so-called `topological optimization' (cf. Novotny and J.~Soko\l owski \cite{NoSo13}). The study of the energy integral also allows to investigate the capacity and then to deduce asymptotic expansions for the eigenvalues of the Dirichlet Laplacian in perforated domains (see, {\it e.g.}, Courtois \cite{Co95} and Abatangelo, Felli, Hillairet, and L\'ena \cite{AbFeHiLe16}). The authors plan to present a detailed analysis on this subject in forthcoming papers.

\vspace{\baselineskip}

\section*{Acknowledgement}

M. Dalla Riva and P. Musolino acknowledge the support of `Progetto di Ateneo: Singular perturbation problems for differential operators -- CPDA120171/12' - University of Padova. M. Dalla Riva also acknowledges the support of HORIZON 2020 MSC EF project FAANon (grant agreement MSCA-IF-2014-EF - 654795) at the University of Aberystwyth, UK. P.~Musolino also acknowledges the support of `INdAM GNAMPA Project 2015 - Un approccio funzionale analitico per problemi di perturbazione singolare e di omogeneizzazione' and of an INdAM Research Fellowship. Part of the work has been carried out while P.~Musolino was visiting the `D\'{e}partement de math\'{e}matiques et applications' of the `\'{E}cole normale sup\'{e}rieure, Paris'. P.~Musolino wishes to thank the `D\'{e}partement de math\'ematiques et applications' and in particular V.~Bonnaillie-No\"el for the kind hospitality. P.~Musolino is a S\^er CYMRU II COFUND fellow, also supported by the `S\^er Cymru National Research Network for Low Carbon, Energy and Environment'.

\end{document}